\documentclass[11pt]{article}

\usepackage{amsfonts,amssymb,amsmath,latexsym,indentfirst,cite}

\newtheorem{rema}{Remark}[section]
\newtheorem{defi}{Definition}[section]
\newtheorem{lemm}{Lemma}[section]
\newtheorem{theo}{Theorem}[section]
\newtheorem{coro}{Corollary}[section]

\newcommand{\Z}[1][]{\ensuremath{{\mathbb{Z}^{#1}} }}

\newcommand{\C}[1][]{\ensuremath{{\mathbb{C}^{#1}} }}
\newcommand{\R}{\ensuremath{{\mathbb{R}} }}

\newcommand{\peq}{\hspace*{0.10in}}
\newcommand{\ppeq}{\hspace*{0.05in}}

\newcommand{\fim}{\rightline{$\blacksquare$}}

\author{Luiz Gustavo Farah$^1$\footnote{Partially supported by FAPESP-Brazil.}
and Ademir Pastor$^2$\footnote{Supported by CNPq-Brazil under grant 152234/2007-1.}.\\
\\
$^1$Department of Mathematics - IMECC/UNICAMP\\
C.P. 6065\\
CEP 13083-970, Campinas, SP, BRAZIL.\\
 and\\
$^2$Instituto de Matem\'atica Pura e Aplicada (IMPA)\\
Estrada Dona Castorina 110, \\
CEP 22460-320, Rio de Janeiro, RJ, BRAZIL. }

\title{On the periodic Schr\"odinger-Boussinesq System }
\date{}

\begin{document}
\maketitle

\begin{abstract}
We study the local and global well-posedness of  the periodic
boundary value problem for the nonlinear Schr\"odinger-Boussinesq
system. The existence of periodic pulses as well as the stability of
such solutions are also considered.
\end{abstract}

\section{Introduction}
In this paper we consider the periodic Schr\"odinger-Boussinesq
system (hereafter referred to as the $SB$-system)
\begin{eqnarray}\label{SB}
\left\{
\begin{array}{l}
iu_{t}+u_{xx}=\alpha vu, \\
v_{tt}-v_{xx}+v_{xxxx}=\beta (|u|^2)_{xx},
\end{array} \right.
\end{eqnarray}
where $t>0$, $x\in [0,L]$, for some $L>0$, and $\alpha,\beta$ are
real constants  .

Here $u$ and $v$ are respectively a complex-valued and a real-valued
function defined in space-time $[0,L] \times \R$. The $SB$-system is
considered as a model of interactions between short and intermediate
long waves, which is derived in describing the dynamics of Langmuir
soliton formation and interaction in a plasma \cite{M} and diatomic
lattice system \cite{YS}. The short wave term $u(x,t): [0,L]\times\R
\rightarrow \C$ is described by a Schr\"odinger type equation with a
potential $v(x,t): [0,L]\times\R \rightarrow \R$ satisfying some
sort of Boussinesq equation and representing the intermediate long
wave.

The nonlinear Schr\"odinger (NLS) equation models a wide range of
physical phenomena including self-focusing of optical beams in
nonlinear media, propagation of Langmuir waves in plasmas, etc. For
an introduction in this topic, we refer the reader to \cite{LP}.
Boussinesq equation as a model of long waves was originally derived
by Boussinesq \cite{BOU} in his study of nonlinear, dispersive wave
propagation. We should remark that it was the first equation
proposed in the literature to describe this kind of physical
phenomena. This equation was also used by Zakharov \cite{Z} as a
model of nonlinear string and by Falk \textit{et al} \cite{FLS} in
their study of shape-memory alloys.

Our first aim here is to study the well-posedness of the periodic
boundary value problem (BVP) for the $SB$-system (\ref{SB}), that
is,  we are interested in the solvability of system \eqref{SB}
subject to the initial conditions
\begin{equation}\label{0SB}
u(x,0)=u_0(x); \ppeq v(x,0)=v_0(x); \ppeq v_t(x,0)=(v_1)_x(x).
\end{equation}

Concerning the local well-posedness question, some results has been
obtained for the $SB$-system (\ref{SB}) in the continuous case.
Linares and Navas \cite{LN} proved that (\ref{SB}) is locally
well-posedness for initial data $u_0 \in L^{2}(\R)$, $v_0 \in
L^{2}(\R)$, $v_1 = h_{x}$ with $h \in {H}^{-1}(\R)$ and $u_0 \in
H^{1}(\R)$, $v_0 \in H^{1}(\R)$, $v_1 = h_{x}$ with $h \in
{L}^{2}(\R)$. Moreover, by using some conservations laws, in the
latter case the solutions can be extended globally. Yongqian
\cite{Y} established a similar result when $u_0 \in H^{s}(\R)$, $v_0
\in H^{s}(\R)$, $v_1= h_{xx}$ with $h\in H^{s}(\R)$ for $s\geq 0$
and assuming $s\geq 1$ these solutions are global. Finally, Farah
\cite{LG2} proved local well-posedness for initial data
$(u_0,v_0,v_1)\in H^k(\R)\times H^s(\R)\times H^{s-1}(\R)$ provided
\begin{enumerate}
\item [$(i)$] $|k|-1/2<s<1/2+2k$ for $k\leq 0$,
\item [$(ii)$] $k-1/2<s< 1/2+k$ for $k> 0$.
\end{enumerate}

In particular, local well-posedness holds for initial data
$(u_0,v_0,v_1)\in H^{s}(\R)\times H^{s}(\R) \times H^{s-1}(\R)$ with
$s>-1/4$. Moreover when $s=0$ the solution is global. We should
mention that, in fact, it is possible to obtain global
well-posedness for $s\geq 0$ in the continuous case. This can be
proved using the arguments introduced by Bourgain \cite{B5} (see
also Angulo \textit{et al.} \cite{AMP1}). In the proof of Theorem
\ref{t1.5} below we also apply these techniques for the periodic
$SB$-system (\ref{SB})-(\ref{0SB}).

The local well-posedness for single dispersive equations with
quadratic nonlinearities has been extensively studied in Sobolev
spaces. The proof of these results are based in the Fourier
restriction norm approach introduced by Bourgain \cite{B} in his
study of the nonlinear Schr\"odinger (NLS) equation
$iu_{t}+u_{xx}+u|u|^{p-2}=0$, with $p\geq 3$  and the Korteweg-de
Vries (KdV) equation $u_{t}+u_{xxx}+u_{x}u=0$. This method was
further developed by Kenig, Ponce and Vega in \cite{KPV1} for the
KdV equation and \cite{KPV2} for the quadratics nonlinear
Schr\"odinger equations
\begin{eqnarray*}
iu_{t}+u_{xx}+F_j(u,\bar{u})=0, \,\,\, j=1,2,3,
\end{eqnarray*}
where $\bar{u}$ denotes the complex conjugate of $u$ and
$F_1(u,\bar{u})=u^2$,  $F_2(u,\bar{u})=u\bar{u}$,
$F_3(u,\bar{u})={\bar{u}}^2$ in one spatial dimension and in
spatially continuous and periodic case.

The original Bourgain method makes extensive use of the Strichartz
inequalities in order to derive the bilinear estimates corresponding
to the nonlinearity. On the other hand, Kenig, Ponce and Vega
simplified Bourgain's proof and improved the bilinear estimates
using only elementary techniques, such as Cauchy-Schwarz inequality
and simple calculus inequalities.

This same kind of technique was used by Farah \cite{LG4} for the
Boussinesq equation. However, we do not have good cancellations on
the Boussinesq symbol. To overcome this difficulty, we observed that
the dispersion in the Boussinesq case is given by the symbol
$\sqrt{{\xi}^2+{\xi}^4}$ and this is, in some sense, related with
the Schr\"odinger symbol (see Lemma \ref{l3.3} below). Therefore, we
can modify the symbols and work only with the algebraic relations
for the Schr\"odinger equation already used in Kenig, Ponce and Vega
\cite{KPV2} in order to derive our relevant bilinear estimates.

To describe our results we define next the $X^S_{s,b}$ and
$X^B_{s,b}$ spaces related respectively to the Schr\"odinger and
Boussinesq equations. For the first equation, this spaces were
introduced in \cite{B}. In the case of Boussinesq equation, the
$X^B_{s,b}$, were first defined by Fang and Grillakis \cite{FG} for
the Boussinesq-type equations in the periodic case. Using these
spaces and following Bourgain's argument introduced in \cite{B} they
proved local well-posedness for the BVP
\begin{eqnarray*}
\left\{
\begin{array}{l}
u_{tt}-u_{xx}+u_{xxxx}+\partial_x^2[f(u)]=0,\\
u_x,0)=u_0(x), \ppeq u_t(x,0)=(u_1)_x(x),
\end{array} \right.
\end{eqnarray*}
where $u_0 \in H_{per}^{s}$, $u_1 \in H_{per}^{-2+s}$, with $0\leq
s\leq 1$ and the nonlinearity $f$ satisfying $|f(u)|\leq c|u|^{p}$,
with $1<p<\frac{3-2s}{1-2s}$ if $0\leq s<\frac{1}{2}$ and
$1<p<\infty$ if $\frac{1}{2}\leq s\leq 1$. Moreover, if $u_0 \in
H_{per}^{1}$, $u_1 \in H_{per}^{-1}$ and  $f(u)=
\lambda|u|^{q-1}u-|u|^{p-1}u$, with $1<q<p$ and $\lambda \in \R$
then the solution is global.

Next we give the precise definition of the $X^S_{s,b}$ and
$X^B_{s,b}$ spaces used in the sequel.
\begin{defi}\label{BOS} Let $\mathcal{Y}$ be the space of
functions $F(\cdot)$ such that
\begin{enumerate}
\item [$(i)$] $F: [0,L]\times \R \rightarrow \C$.
\item [$(ii)$] $F(x,\cdot)\in S(\R)$ for each $x\in [0,L]$.
\item [$(iii)$] $F(\cdot,t)\in C^{\infty}([0,L])$ for each $t\in \R$.
\end{enumerate}
For $s,b \in \R$, $X^S_{s,b}$ and $X^B_{s,b}$ denotes,
respectively, the completion of $\mathcal{Y}$ with respect to the norm
\begin{eqnarray}
\|F\|_{X^S_{s,b}}&=&\|\langle\tau+(2\pi n/L)^2\rangle^b\langle n\rangle^s \widetilde{F}\|_{l^{2}_{n}L^{2}_{\tau}},\\
\|F\|_{X^B_{s,b}}&=&\|\langle|\tau|-\gamma_L(n)\rangle^b\langle n\rangle^s \widetilde{F}\|_{l^{2}_{n}L^{2}_{\tau}},
\end{eqnarray}
where $\sim$ denotes the time-space Fourier transform, $\langle
a\rangle\equiv 1+|a|$ and $\gamma_L(n)\equiv(2\pi/L)^2\sqrt{{n}^2+{n}^4}$.
\end{defi}

We will also need the localized $X_{s,b}$ spaces defined as follows:

\begin{defi}
Let $I$ be a time interval. For $s,b \in \R$, $X_{s,b}^{S,I}$ and $X^{B,I}_{s,b}$ denotes the space endowed with the norm
\begin{eqnarray*}
\|u\|_{X_{s,b}^{S,I}}&=&\inf_{w\in X^{S}_{s,b}}\left\{\|w\|_{X^{S}_{s,b}}:w(t)=u(t) \textrm{ on }  I\right\},\\
\|u\|_{X_{s,b}^{B,I}}&=&\inf_{w\in X^{B}_{s,b}}\left\{\|w\|_{X^{B}_{s,b}}:w(t)=u(t) \textrm{ on }  I\right\}.
\end{eqnarray*}
\end{defi}

Now we state our main results concerning well-posedness.

\begin{theo}\label{t1.1}
Let $s\geq 0$ and $1/4<a<1/2<b$. Then, there exists $c>0$,
depending only on $a,b,s$, such that
\begin{enumerate}
 \item [$(i)$]  \peq
$\left\|
uv
\right\|_{X^S_{s,-a}}\leq c\left\|u\right\|_{X^S_{s,b}}\left\|v\right\|_{X^B_{s,b}}.
$
\item [$(ii)$]
\peq
$\left\|
u_1\bar{u}_2
\right\|_{X^B_{s,-a}}\leq c\left\|u_1\right\|_{X^S_{s,b}}\left\|u_2\right\|_{X^S_{s,b}}.
$
\end{enumerate}
\end{theo}

\begin{theo}\label{t1.2}
Let $s\geq 0$. Then for any $(u_0,v_0,v_1)\in H_{per}^s([0,L])\times
H_{per}^s([0,L])\times H_{per}^{s-1}([0,L])$ there exist
$T=T(\|u_0\|_{H_{per}^s},\|v_0\|_{H_{per}^s},\|v_1\|_{H_{per}^{s-1}})$,
$b>1/2$ and a unique solution $(u,v)$ of the BVP
(\ref{SB})--(\ref{0SB}), satisfying
\begin{eqnarray*}
u\in C([0,T]:H_{per}^s([0,L]))\cap X_{s,b}^{S,[0,T]}\textrm{ and }v\in C([0,T]:H_{per}^s([0,L]))\cap X_{s,b}^{B,[0,T]}.
\end{eqnarray*}
Moreover, the map $(u_0,v_0,v_1)\mapsto (u(t),v(t))$ is locally Lipschitz from
$H_{per}^s([0,L])\times H_{per}^s([0,L])\times H_{per}^{s-1}([0,L])$ into $C([0,T]:H_{per}^s([0,L])\times H_{per}^s([0,L]))$.
\end{theo}

We also obtain counter-examples for the bilinear estimates stated in
Theorem \ref{t1.1}.
\begin{theo}\label{t1.3}
\end{theo}
\vspace*{-0.25in}
\emph{\begin{enumerate}
 \item [$(i)$]  The estimate
\begin{equation}\label{CE1}
\left\|uv\right\|_{X^S_{k,-a}}\leq c\left\|u\right\|_{X^S_{k,b}}\left\|v\right\|_{X^B_{s,b}}
\end{equation}
holds only if $k\leq s$.
\item [$(ii)$]  The estimate
$$\left\|uv\right\|_{X^S_{k,-a}}\leq c\left\|u\right\|_{X^S_{k,b}}\left\|v\right\|_{X^B_{s,b}}
$$ holds only if $k+s\geq 0$.
\item [$(iii)$] The estimate
\begin{equation}\label{CE2}
\left\|u_1\bar{u}_2\right\|_{X^B_{s,-a}}\leq c\left\|u_1\right\|_{X^S_{k,b}}\left\|u_2\right\|_{X^S_{k,b}}
\end{equation}
holds only if $s\leq k$.
\end{enumerate}
}

Theorem \ref{t1.3} has an important consequence. It shows that our
local well-posed result is sharp, in the sense that it cannot be
improved using the spaces $X^S_{s,b}$ and $X^B_{s,b}$. This
situation is very different from the continuous case obtained in
Farah \cite{LG2} where we have local well-posedness for initial data
in different Sobolev spaces with negative indices.\\

Next we obtain bilinear estimates for the case $s=0$ and
$b,b_1<1/2$. These estimates will be useful to establish the
existence of global solutions.

\begin{theo}\label{t1.4}
Let $a,a_1,b,b_1>1/4$, then there exists $c>0$ depending only on $a,a_1,b,b_1$ such that
\begin{enumerate}
 \item [(i)]  \peq
$\left\|
uv
\right\|_{X^S_{0,-a_1}}\leq c\left\|u\right\|_{X^S_{0,b_1}}\left\|v\right\|_{X^B_{0,b}}.
$
\item [(ii)]
\peq
$\left\|
u_1\bar{u}_2
\right\|_{X^B_{0,-a}}\leq c\left\|u_1\right\|_{X^S_{0,b_1}}\left\|u_2\right\|_{X^S_{0,b_1}}.
$
\end{enumerate}
\end{theo}

The bilinear estimates in Theorem \ref{t1.4} are the essential tools
to prove the global result. It asserts that the local solution given
by Theorem \ref{t1.2} is in fact a global one, for all $s\geq0$.

\begin{theo}\label{t1.5}
Let $s\geq 0 $. Then, the BVP \eqref{SB}--\eqref{0SB} is globally
well-posed for data $(u_0,v_0,v_1)\in H_{per}^{s}([0,L])\times
H_{per}^{s}([0,L]) \times H_{per}^{s-1}([0,L])$. Moreover, the
solution $(u,v)$ satisfies, for all $t>0$,
\begin{eqnarray*}
 \|v(t)\|_{H_{per}^{s}}+ \|(-\Delta)^{-1/2}v_t(t)\|_{H_{per}^{s-1}}\lesssim e^{((\ln{2})\|u_0\|_{H_{per}^{s}}^2t)} \max{\{\|v_0,v_1\|_{\mathfrak{B^s}},\|u_0\|_{H_{per}^{s}}\}}.
\end{eqnarray*}
\end{theo}

The argument used to prove Theorem \ref{t1.5} follows the ideas
introduced by Colliander, Holmer, Tzirakis \cite{CHT} to deal with
the Zakharov system. The intuition for this Theorem comes from the
fact that the nonlinearity for the second equation of the
$SB$-system (\ref{SB}) depends only on the first equation.
Therefore, noting that the bilinear estimates given in Theorem
\ref{t1.4} hold for $a,a_1,b,b_1<1/2$, it is possible to show that
the time existence depends only on the $\|u_0\|_{L_{per}^2}$. But
since this norm is conserved by the flow, we obtain a global
solution.

Our second aim is to study existence and orbital (nonlinear)
stability of periodic traveling waves. These two questions  are very
important in the understanding of the dynamic of the system under
consideration.

The stability study of traveling waves has been extensively studied
for the whole Euclidean space case (solitary waves), whereas the
study under periodic boundary conditions has been started quite
recently and few works are available in the current literature. To
cite a few important contributions, in \cite{Angulo} Angulo studied
the orbital stability of {\it dnoidal} wave solutions for the cubic
Schr\"odinger and modified Korteweg-de Vries equations; his method
of proofs follows the pioneers ideas of Benjamin, Bona and
Weinstein. In \cite{AnBoSci}, Angulo \emph{et al.} gave a complete
stability study of {\it cnoidal} wave solutions for the Korteweg-de
Vries equation, adapting to the the periodic context the abstract
theory developed in \cite{GSS}. For others well-known equations and
systems see e.g. \cite{AL}, \cite{AMP1}, \cite{deka}, \cite{hailki},
\cite{NP1} (and references therein).

One of the main reasons why the stability study in the periodic case
has been received little attention, lies on the needed spectral
theory associated with some linearized operator. Indeed, to fix
ideas, suppose we have a Schr\"odinger type operator
$\mathcal{L}=-\frac{d^2}{dx^2}+q(x)$, where $q(x)$ is a smooth real
potential. If $q(x)$ and $\phi$ are rapidly decaying  smooth
functions such that $\mathcal{L}\phi=0$ and assuming that $\phi$ has
exactly two zeros on the whole real line, then it follows
immediately from Sturn-Liouville's theory that zero is the third
eigenvalue of  operator $\mathcal{L}$ and it is a simple eigenvalue.
On the other hand, if $q(x)$ is a periodic function with period
$L>0$ and $\phi$ is also $L$-periodic such that $\mathcal{L}\phi=0$
and has exactly two zeros on the interval $[0,L)$ then from
Floquet's theory, the eigenvalue zero is the second or the third one
(see e.g. \cite{EASTHAM}). In most cases, it is a hard task to
decide when zero is the second or the third eigenvalue. As a
consequence, most of the current papers deal with explicit
solutions. This is the case of the present paper.

In general, the studied dispersive equations admits periodic
explicit solutions depending on the Jacobian elliptic functions
(dnoidal, cnoidal and snoidal type). So, the main idea to obtain the
spectral properties for the linearized operator is to reduce matter
to some known periodic eigenvalue problem. The most popular one
deals with the periodic eigenvalue problem associated with the
Lam\'e operator
\begin{equation}  \label{lame}
\mathcal{L}_{Lame}:=-\frac{d^2}{dx^2}+n(n+1)sn^2(x;k),
\end{equation}
for some determined value of $n\in\mathbb{N}$ (see e.g.
\cite{Angulo}, \cite{AnBoSci}, \cite{AL}, \cite{NP1}).

Here, we will consider $\alpha=\beta=-1$ in \eqref{SB} and look for
solutions of the form
\begin{equation}\label{s22}
u(x,t)=e^{i\omega t} \psi_{\omega}(x), \qquad v(x,t)=\phi_\omega(x),
\end{equation}
where $\omega$ is a real parameter and $\psi_\omega,
\phi_\omega:\mathbb{R}\rightarrow\mathbb{R}$ are $L$-periodic
functions with a period $L>0$. Then, substituting this waveform into
the system and integrating twice the second equation in the obtained
system, we have
\begin{eqnarray}\label{s44}
\left\{
\begin{array}{l}
\psi_\omega''-\omega\psi_\omega+\psi_\omega\phi_\omega=0, \\
\phi_\omega''-\phi_\omega+\psi_\omega^2=0.\\
\end{array} \right.
\end{eqnarray}
To reduce system \eqref{s44} to a single ordinary differential
equation, we assume $\omega=1$ and $\psi_\omega=\phi_\omega=\psi$,
so that it reduces to
\begin{equation}\label{s55}
\psi''-\psi+\psi^2=0.
\end{equation}

Before proceeding, we point out that existence and stability of
hyperbolic-secant-type solitary waves were recently considered in
\cite{hakkaev}. The author has proved a orbital stability result by
using the abstract theory contained  in \cite{GSS}, taking the
advantage of the spectral properties established in \cite{Lopes}.

In the periodic approach, it is not difficult to prove that
\eqref{s55} has a periodic solution of {\it cnoidal} type, namely,
\begin{equation}   \label{cnwave}
\displaystyle\psi(x)=
\beta_2+(\beta_3-\beta_2)cn^2\left(\sqrt{\frac{\beta_3-\beta_1}{6}}x;k\right),
\qquad k^2=\frac{\beta_3-\beta_2}{\beta_3-\beta_1}
\end{equation}
where $cn(\cdot,k)$ denotes the Jacobian elliptic function of
dnoidal type and $\beta_1,\beta_2,\beta_3$ are real parameters.

Our main theorem concerned with the orbital stability of cnoidal
waves reads as follows:

\begin{theo}  \label{stability}
Let $\psi$ be the cnoidal wave solution given in \eqref{cnwave}.
Then, the periodic traveling wave $(e^{it}\psi,\psi)$ is orbitally
stable in the energy space $X= H_{per}^1([0,L])\times
H_{per}^1([0,L])\times L_{per}^2([0,L])$ by the flow of system
\eqref{SB}.
\end{theo}

To prove Theorem \ref{stability}, we shall employ the classical
theory developed by Grillakis, Shatah and Strauss \cite{GSS}. To do
so, we first observe that system \eqref{SB} (with $\alpha=\beta=-1$)
can be written in Hamiltonian form (see \eqref{sa1}). We point out
that although the operator $J$ in \eqref{sa2} is not onto,  along
the lines of proofs in \cite{GSS} the stability result still holds
(see also \cite{hakkaev}, \cite{strauss}).

Our strategy to get the needed spectral properties is to combine the
results in \cite{AL}, which are essentially proved  from well-known
results for the Lam\'e operator in \eqref{lame}, with the min-max
principle for the eigenvalues characterization.

Finally, we also obtain periodic traveling waves  for $\omega\neq1$.
Our idea  is simple: once obtained the cnoidal solution for
$\omega=1$, we employ the Implicit Function Theorem combined with
spectral properties related with the  linearized operator  to extend
our range of parameters for $\omega$ near 1.

The plan of this paper is as follows: in Section 2, we introduce
some  notation and state important propositions that we will use
throughout the paper. The proof of the bilinear estimates and the
relevant counter examples are given in Sections 3 and 4,
respectively. In Section 5 we prove Theorem \ref{t1.5}. Finally, the
stability questions are treated in Sections 6.

\section{Notations and Preliminaries}

In what follows we use $a \lesssim b$ to say that $a \leq C b$ for
some constant $C>0$. Also, we denote $a \sim b$ when, $a \lesssim b$
and $b \lesssim a$. We write $a\ll b$ to denote an estimate of the
form $a\leq c b$ for some small constant $c>0$. In addition, $a+$
means that there exists $\varepsilon>0$ such that
$a+=a+\varepsilon$.

Let us recall some properties of $L$-periodic functions. For a
detailed presentation of the spaces of periodic functions and its
properties we refer the reader, for instance, to \cite{IoIo}. We
define the Fourier transform of $f\in L^1([0,L])$ by
\begin{equation*}
 \hat{f}(n)=\dfrac{1}{L}\int_{0}^{L}e^{-2\pi i\frac{x}{L}n}f(x)dx.
\end{equation*}

For $f$ in an appropriate class of functions we have
$f=(\widehat{f})^{\vee}$, where for a sequence
$s=\{s_n\}_{n\in\mathbb{Z}}$, the symbol ${}^\vee$ denotes the
inverse Fourier transform of $s$ given by
\begin{equation*}
 (s)^{\vee}(x)=\sum_{n\in \mathbb{Z}}e^{2\pi i\frac{x}{L}n}s_n
\end{equation*}

Moreover, we have the Plancherel identity
\begin{equation*}
 \|f\|_{L^2_{per}}=\|\hat{f}\|_{l^2_{n}}.
\end{equation*}

The periodic Sobolev space $H_{per}^s([0,L])$ is defined to be space
of all periodic distributions such that
\begin{equation*}
\|f\|_{H^s_{per}}:=\|\langle
n\rangle^s\hat{f}(n)\|_{l^2_{n}}<\infty.
\end{equation*}

Moreover, the operator $(-\Delta)^{-1/2}$ is defined, via Fourier transform,
by
$$
[(-\Delta)^{-1/2}f]^{\wedge}(n)=|n|^{-1}\hat{f}(n)\qquad n\neq0.
$$

Next, we recall some facts on the linear Schr\"odinger and
Boussinesq equations. Consider the free Schr\"odinger equation
\begin{equation}\label{LNLS}
iu_{t}+u_{xx}=0.
\end{equation}
It is easy to see that the solution of \eqref{LNLS}, with initial
data $u(0)=u_0$, is given by the formula
\begin{equation}\label{GUS}
u(t)=U(t)u_0,
\end{equation}
where
\begin{eqnarray*}
U(t)u_0&=&\left( e^{-(2\pi/L)^2 itn^2}\widehat{u}_0(n)\right)^{\vee}.\\
\end{eqnarray*}

On the other hand, for the linear Boussinesq equation
\begin{equation}\label{LB}
v_{tt}-v_{xx}+v_{xxxx}=0
\end{equation}
it is well-known that the solution, with initial data $v(0)=v_0$ and
$v_t(0)=(v_1)_x$, is given by
\begin{equation}\label{GUB}
u(t)=V_c(t)v_0+V_s(t)(v_1)_x,
\end{equation}
where
\begin{eqnarray*}
V_c(t)v_0&=&\left( \frac{e^{(2\pi/L)^2 it\sqrt{{n}^2+{n}^4}}+ e^{-(2\pi/L)^2 it\sqrt{{n}^2+{n}^4}}}{2}\widehat{v_0}(n)\right)^{\vee}\\
V_s(t)(v_1)_x&=&\left( \frac{e^{(2\pi/L)^2 it\sqrt{{n}^2+{n}^4}}- e^{-(2\pi/L)^2 it\sqrt{{n}^2+{n}^4}}}{2i\sqrt{{n}^2+{n}^4}}\widehat{(v_1)_x}(n)\right)^{\vee}.
\end{eqnarray*}

As a consequence, by Duhamel's Principle the solution of
\eqref{SB}--\eqref{0SB}, is equivalent to
\begin{eqnarray}\label{INT}
\begin{split}
u(t)=& U(t)u_0-i\int_{0}^{t}U(t-t')(\alpha vu)(t')dt'\\
v(t)=&
V_c(t)v_0+V_s(t)(v_1)_x+\int_{0}^{t}V_s(t-t')(\beta|u|^2)_{xx}(t')dt'.
\end{split}
\end{eqnarray}

Let $\theta$ be a cutoff function satisfying $\theta \in
C^{\infty}_{0}(\R)$, $0\leq \theta \leq 1$, $\theta \equiv 1$ in
$[-1,1]$, supp$(\theta) \subseteq [-2,2]$ and for $0<T\leq 1$ define
$\theta_T(t)=\theta(t/T)$. In fact, to work in the $X^S_{s,b}$ and
$X^B_{s,b}$ we consider another versions of (\ref{INT}), that is
\begin{eqnarray}\label{INT2}
\begin{split}
u(t)=& \theta(t)U(t)u_0-i\theta_T(t)\int_{0}^{t}U(t-t')(\alpha vu)(t')dt'\\
v(t)=& \theta(t)\left(V_c(t)v_0+V_s(t)(v_1)_x\right)+\theta_T(t)
\int_{0}^{t}V_s(t-t')(\beta|u|^2)_{xx}(t')dt'.
\end{split}
\end{eqnarray}
and
\begin{eqnarray}\label{INT4}
\begin{split}
u(t)=& \theta_T(t)U(t)u_0-i\theta_T(t)\int_{0}^{t}U(t-t')(\alpha vu)(t')dt'\\
v(t)=& \theta_T(t)\left(V_c(t)v_0+V_s(t)(v_1)_x\right)+\theta_T(t)
\int_{0}^{t}V_s(t-t')(\beta|u|^2)_{xx}(t')dt'.
\end{split}
\end{eqnarray}

We will use equation (\ref{INT2}) (resp. (\ref{INT4})) to study the
local (resp. global) well-posedness problem associated to
(\ref{SB})--\eqref{0SB}.

Note that the integral equations (\ref{INT2}) and (\ref{INT4}) are
defined for all $(t,x)\in \R^2$. Moreover, if $(u,v)$ is a solution
of (\ref{INT2}) or (\ref{INT4}) then
$(\tilde{u},\tilde{v})=(u|_{[0,T]}, v|_{[0,T]})$ will be a solution
of (\ref{INT}) in $[0,T]$.

Before proceeding to the group and integral estimates for (\ref{INT2}) and (\ref{INT4}) we introduce the norm
\begin{equation*}
\|v_0,v_1\|_{\mathfrak{B}^s}^2\equiv \|v_0\|_{H_{per}^s([0,L])}^2+ \|v_1\|_{H_{per}^{s-1}([0,L])}^2.
\end{equation*}

For simplicity we denote $\mathfrak{B}^0$ by $\mathfrak{B}$ and, for functions of $t$, we use the shorthand
\begin{equation*}
\|v(t)\|_{\mathfrak{B}^s}^2\equiv \|v(t)\|_{H_{per}^s([0,L])}^2+ \|(-\Delta)^{-1/2}v_t(t)\|_{H_{per}^{s-1}([0,L])}^2.
\end{equation*}

The following three lemmas are standard in this context. Although we are studying the periodic case, the proofs are essentially the same of the continuous setting. We refer the reader to Farah \cite{LG2} for the details.
\begin{lemm}[Group estimates]\label{l21}
Let $L=2\pi$ and $0<T\leq 1$.
\begin{enumerate}
\item [(a)]Linear Schr\"odinger equation
\begin{itemize}
\item [(i)]$\|U(t)u_0\|_{C(\R:H_{per}^s)}=\|u_0\|_{H_{per}^s}.$
\item [(ii)] If $0\leq b_1\leq 1$, then
\begin{equation*}
\|\theta_T(t)U(t)u_0\|_{X^S_{s,b_1}}\lesssim T^{1/2-b_1}\|u_0\|_{H_{per}^s}.
\end{equation*}
\end{itemize}
\item [(b)]Linear Boussinesq equation
\begin{itemize}
\item [(i)]$\|V_c(t)v_0+V_s(t)(v_1)_x\|_{C(\R:H_{per}^s)}\leq \|v_0\|_{H_{per}^s}+\|v_1\|_{H_{per}^{s-1}}.$
\item [(ii)]$\|V_c(t)v_0+V_s(t)(v_1)_x\|_{C(\R:\mathfrak{B})}= \|v_0,v_1\|_{\mathfrak{B}}.$
\item [(iii)]If $0\leq b\leq 1$, then
\begin{equation*}
\|\theta_T(t)\left(V_c(t)v_0+V_s(t)(v_1)_x\right)\|_{X^B_{s,b}}\lesssim T^{1/2-b}\left(\|v_0\|_{H_{per}^s}+\|v_1\|_{H_{per}^{s-1}}\right).
\end{equation*}
\end{itemize}
\end{enumerate}
\end{lemm}

Next we estimate the integral parts of (\ref{INT2}).

\begin{lemm}[Integral estimates]\label{l22}
Let $L=2\pi$ and $0<T\leq 1$.
\begin{enumerate}
\item [(a)] Nonhomogeneous linear Schr\"odinger equation
\begin{itemize}
\item [(i)]If $0\leq a_1 <1/2$ then
\begin{equation*}
\left\|\int_{0}^{t}U(t-t')z(t')dt'\right\|_{C([0,T]:H_{per}^s)}\lesssim T^{1/2-a_1}\|z\|_{X^S_{s,-a_1}}.
\end{equation*}
\item [(ii)]If $0\leq a_1 <1/2$, $0\leq b_1$ and $a_1+b_1\leq 1$ then
\begin{equation*}
\left\|\theta_T(t)\int_{0}^{t}U(t-t')z(t')dt'\right\|_{X^S_{s,b_1}}\lesssim T^{1-a_1-b_1}\|z\|_{X^S_{s,-a_1}}.
\end{equation*}
\end{itemize}
\item [(b)]Nonhomogeneous linear Boussinesq equation
\begin{itemize}
\item [(i)]If $0\leq a <1/2$ then
\begin{equation*}
\left\|\int_{0}^{t}V_s(t-t')z_{xx}(t')dt'\right\|_{C([0,T]:\mathfrak{B}^s)}\lesssim T^{1/2-a}\|z\|_{X^B_{s,-a}}.
\end{equation*}
\item [(ii)]If $0\leq a <1/2$, $0\leq b$ and $a+b\leq 1$ then
\begin{equation*}
\left\|\theta_T(t)\int_{0}^{t}V_s(t-t')z_{xx}(t')dt'\right\|_{X^B_{s,b}}\lesssim T^{1-a-b}\|z\|_{X^B_{s,-a}}.
\end{equation*}
\end{itemize}
\end{enumerate}
\end{lemm}

We also know the following embeddeding concerning the $X^S_{s,b}$ and $X^B_{s,b}$ spaces.

\begin{lemm}\label{l23}
Let $b>\frac{1}{2}$. There exists $c>0$, depending only on $b$, such that
\begin{eqnarray*}
\|u\|_{C(\R:H_{per}^s)}&\leq& c\|u\|_{X^B_{s,b}}\\
\|u\|_{C(\R:H_{per}^s)}&\leq& c\|u\|_{X^S_{s,b}}.
\end{eqnarray*}
\end{lemm}

We finish this section with the following standard Bourgain-Strichartz estimates.
\begin{lemm}\label{l24}
Let $u\in L^3_{x,t}$, therefore
\begin{equation*}
\|u\|_{L^3_{x,t}}\leq c\min\{\|u\|_{X^S_{0,1/4+}},\|u\|_{X^B_{0,1/4+}}\}.
\end{equation*}
\end{lemm}

\textbf{Proof. }This estimate is easily obtained by interpolating between
\begin{itemize}
\item   $\|u\|_{L^4_{x,t}}\leq c\min\{\|u\|_{X^S_{0,3/8+}},\|u\|_{X^B_{0,3/8+}}\}$ (See Bougain \cite{B} and Fang and Grillakis \cite{FG}).
\item   $\|u\|_{L^2_{x,t}}= \|u\|_{X^S_{0,0}}=\|u\|_{X^B_{0,0}}$ (by definition).
\end{itemize}
\fim

\begin{rema}
To simplify our  well-posedness analysis we will assume $L=2\pi$. We will return to an arbitrarily $L>0$ in Section 6, where we study stability questions.
\end{rema}


\section{Bilinear estimates}

First we state some elementary calculus inequalities that will be useful later.
\begin{lemm}\label{l3.1}
For $p, q>0$ and $r=\min\{ p, q, p+q-1\}$ with $p+q>1$, there exists $c>0$ such that
\begin{equation}\label{CI1}
\int_{-\infty}^{+\infty}\dfrac{dx} {\langle x-\alpha\rangle^{p}\langle x-\beta\rangle^{q}}\leq\dfrac{c} {\langle \alpha-\beta\rangle^{r}}.
\end{equation}
\end{lemm}
\textbf{Proof. } See Lemma 4.2 in \cite{GTV}.\\
\fim

\begin{lemm}\label{l3.2}
If $\gamma>1/2$, then
\begin{equation}\label{CI2}
\sup_{(n,\tau)\in \Z\times\R}\sum_{n_1\in\Z}\dfrac{1}{(1+|\tau\pm n_1(n-n_1)|)^{\gamma}}<\infty.
\end{equation}
\end{lemm}
\textbf{Proof. } See Lemma 5.3 in \cite{KPV2}.\\
\fim

\begin{lemm}\label{l3.3}
There exists $c>0$ such that
\begin{equation}\label{LN}
\dfrac{1}{c}\leq\sup_{x,y\geq 0}\dfrac{1+|x-y|}{1+|x-\sqrt{y^2+y}|}\leq c.
\end{equation}
\end{lemm}

\textbf{Proof. } Since $y\leq\sqrt{y^2+y}\leq y+1/2$ for all $y\geq 0$ a simple computation shows the desired inequalities.\\
\fim
\begin{rema}\label{R1}
In view of the previous lemma we have an equivalent way to estimate the $X_{s,b}^B$-norm, that is
\begin{equation*}
\|u\|_{X_{s,b}^B}\sim\|\langle|\tau|-n^2\rangle^b\langle n\rangle^{s} \widetilde{u}(\tau,n)\|_{l^{2}_{n}L^2_{\tau}}.
\end{equation*}
This equivalence will be important in the proof of Theorem \ref{t1.1}. As we said in the introduction, the Boussinesq symbol $\sqrt{{n}^2+{n}^4}$ does not have good cancellations to make use of Lemma \ref{l3.1}. Therefore, we  modify the symbols as above and work only with the algebraic relations for the Schr\"odinger equation.
\end{rema}

Now we are in position to prove the bilinear estimates stated in Theorem \ref{t1.1}.\\

\textbf{Proof of Theorem \ref{t1.1}}
\begin{enumerate}
\item [$(i)$] For $u\in X^S_{s,b}$ and $v\in X^B_{s,b}$ we define $f(\tau,n)\equiv \langle\tau+n^2\rangle^b\langle n\rangle^{s} \widetilde{u}(\tau,n)$ and $g(\tau,n)\equiv \langle|\tau|-\gamma(n)\rangle^b\langle n\rangle^{s} \widetilde{v}(\tau,n)$. By duality the desired inequality is equivalent to
\begin{equation}\label{DUA}
\left|W(f,g,\phi)\right|\leq c\|f\|_{l^{2}_{n}L^2_{\tau}} \|g\|_{l^{2}_{n}L^2_{\tau}}\|\phi\|_{l^{2}_{n}L^2_{\tau}}
\end{equation}
where
\begin{equation*}
W(f,g,\phi)= \sum_{n,n_1} \int_{\R^2} \dfrac{\langle n\rangle^{s}}{\langle n_1\rangle^{s} \langle n_2\rangle^{s}} \dfrac{g(\tau_1,n_1)f(\tau_2,n_2) \bar{\phi}(\tau,n)}{\langle\sigma\rangle^{a} \langle\sigma_1\rangle^{b} \langle\sigma_2\rangle^{b}} d\tau d\tau_1
\end{equation*}
and\\
\begin{equation}\label{TAU2}
n_2=n-n_1, \peq \tau_2=\tau-\tau_1,
\end{equation}
\centerline{$
\sigma=\tau+n^2,\peq \sigma_1=|\tau_1|-\gamma(n_1),\peq \sigma_2=\tau_2+n_2^2.
$}\\

In view of Remark \ref{R1}, we know that $\langle|\tau_1|-\gamma(n_1)\rangle \thicksim \langle|\tau_1|-n_1^2\rangle$. Therefore splitting the domain of integration into the regions $\{(n, \tau, n_1, \tau_1)\in \R^4: \tau_1<0\}$ and $\{(n, \tau, n_1, \tau_1)\in \R^4: \tau_1\geq 0\}$, it is sufficient to prove inequality (\ref{DUA}) with $W_1(f,g,\phi)$ and $W_2(f,g,\phi)$ instead of $W(f,g,\phi)$, where
 \begin{equation*}
W_1(f,g,\phi)=  \sum_{n,n_1} \int_{\R^2} \dfrac{\langle n\rangle^{s}}{\langle n_1\rangle^{s} \langle n_2\rangle^{s}} \dfrac{g(\tau_1,n_1)f(\tau_2,n_2) \bar{\phi}(\tau,n)}{\langle\sigma\rangle^{a} \langle\tau_1+n_1^2\rangle^{b} \langle\sigma_2\rangle^{b}}d\tau d\tau_1
\end{equation*}
and
\begin{equation*}
W_2(f,g,\phi)=  \sum_{n,n_1} \int_{\R^2} \dfrac{\langle n\rangle^{s}}{\langle n_1\rangle^{s} \langle n_2\rangle^{s}} \dfrac{g(\tau_1,n_1)f(\tau_2,n_2) \bar{\phi}(\tau,n)}{\langle\sigma\rangle^{a} \langle\tau_1-n_1^2\rangle^{b} \langle\sigma_2\rangle^{b}}d\tau d\tau_1.
\end{equation*}

Applying Cauchy-Schwarz and H\"older inequalities it is easy to
see that
\begin{eqnarray*}
|W_1|^2&\leq& \|f\|_{l^{2}_{n}L^2_{\tau}}^2 \|g\|_{l^{2}_{n}L^2_{\tau}}^2 \|\phi\|_{l^{2}_{n}L^2_{\tau}}^2\\
&&\times\left\|\dfrac{\langle n\rangle^{2s}}{\langle\sigma\rangle^{2a}} \sum_{n_1}\int\dfrac{d\tau_1}{\langle n_1\rangle^{2s} \langle n_2\rangle^{2s} \langle\tau_1+ n_1^2\rangle^{2b} \langle\sigma_2\rangle^{2b}}\right\|_{l^{\infty}_{n}L^{\infty}_{\tau}}.
\end{eqnarray*}

Noting that $s\geq 0$ we have
\begin{equation}\label{Xi}
\dfrac{\langle n\rangle^{2s}}{\langle n_1\rangle^{2s} \langle n_2\rangle^{2s}}\leq 1.
\end{equation}

Therefore in view of Lemma \ref{l3.1} it suffices to get bounds for
\begin{eqnarray*}
\sup_{n,\tau}\dfrac{1}{\langle\sigma\rangle^{2a}} \sum_{n_1}\dfrac{1} {\langle\tau+n^2+2n_1^2-2nn_1\rangle^{2b}}.
\end{eqnarray*}
By Lemma \ref{l3.2} this expression is bounded provides $a\geq 0$ and $b>1/4$.

Now we turn to the proof of inequality (\ref{DUA}) with
$W_2(f,g,\phi)$. Using the Cauchy-Schwarz and H\"older inequalities
and duality it is easy to see that
\begin{eqnarray*}
|W_2|^2&\leq&\|f\|_{l^{2}_{n}L^2_{\tau}}^2 \|g\|_{l^{2}_{n}L^2_{\tau}}^2 \|\phi\|_{l^{2}_{n}L^2_{\tau}}^2\\
&&\times\left\|\dfrac{1}{\langle n_2\rangle^{2s}\langle\sigma_2\rangle^{2b}} \sum_{n_1}\int\dfrac{\langle n_1+n_2\rangle^{2s} d\tau_1}{\langle n_1\rangle^{2s}\langle\tau_1-n_1^2\rangle^{2b}   \langle\sigma\rangle^{2a}}\right\|_{l^{\infty}_{n_2}L^{\infty}_{\tau_2}}.
\end{eqnarray*}

Therefore in view of Lemma \ref{l3.1} and (\ref{Xi}) it suffices to get bounds for
\begin{eqnarray*}
\sup_{n_2,\tau_2}\dfrac{1}{\langle\sigma_2\rangle^{2b}} \sum_{n_1}\dfrac{1} {\langle\tau_2+n_2^2+2n_1^2+2n_1n_2\rangle^{2a}}.
\end{eqnarray*}

By Lemma \ref{l3.2} this expression is bounded provides $b\geq 0$ and $a>1/4$.

\item [$(ii)$]For $u_1\in X^S_{s,b}$ and $u_2\in X^S_{s,b}$ we define $f(\tau,n)\equiv \langle\tau+n^2\rangle^b\langle n\rangle^{s} \widetilde{u}_1(\tau,n)$ and $g(\tau,n)\equiv \langle\tau+n^2\rangle^b\langle n\rangle^{s} \widetilde{u}_2(\tau,n)$. By duality the desired inequality is equivalent to
\begin{equation}\label{DUA2}
\left|Z(f,g,\phi)\right|\leq c\|f\|_{l^{2}_{n}L^2_{\tau}} \|g\|_{l^{2}_{n}L^2_{\tau}} \|\phi\|_{l^{2}_{n}L^2_{\tau}}
\end{equation}
where
\begin{equation*}
Z(f,g,\phi)=  \sum_{n,n_1} \int_{\R^2} \dfrac{\langle n\rangle^{s}}{\langle n_1\rangle^{s} \langle n_2\rangle^{s}} \dfrac{h(\tau_1,n_1)f(\tau_2,n_2) \bar{\phi}(\tau,n)}{\langle\sigma\rangle^{a} \langle\sigma_1\rangle^{b} \langle\sigma_2\rangle^{b}}d\tau d\tau_1
\end{equation*}
and\\
\begin{equation*}
h(\tau_1,n_1)=\bar{g}(-\tau_1,-n_1), \peq n_2=n-n_1, \peq \tau_2=\tau-\tau_1,
\end{equation*}
\centerline{$
\sigma=|\tau|-\gamma(n),\peq \sigma_1=\tau_1-n_1^2,\peq \sigma_2=\tau_2+n_2^2.
$}

Therefore applying Lemma \ref{l3.3} and splitting the domain of integration according to the sign of $\tau$ it is sufficient to prove inequality (\ref{DUA2}) with $Z_1(f,g,\phi)$ and $Z_2(f,g,\phi)$ instead of $Z(f,g,\phi)$, where
 \begin{equation*}
Z_1(f,g,\phi)=  \sum_{n,n_1} \int_{\R^2} \dfrac{\langle n\rangle^{s}}{\langle n_1\rangle^{s} \langle n_2\rangle^{s}} \dfrac{h(\tau_1,n_1)f(\tau_2,n_2) \bar{\phi}(\tau,n)}{\langle\tau+n^2\rangle^{a} \langle\sigma_1\rangle^{b} \langle\sigma_2\rangle^{b}}d\tau d\tau_1
\end{equation*}
and
\begin{equation*}
Z_2(f,g,\phi)=  \sum_{n,n_1} \int_{\R^2} \dfrac{\langle n \rangle^{s}}{\langle n_1\rangle^{s} \langle n_2\rangle^{s}} \dfrac{h(\tau_1,n_1)f(\tau_2,n_2) \bar{\phi}(\tau,n)}{\langle\tau-n^2\rangle^{a} \langle\sigma_1\rangle^{b} \langle\sigma_2\rangle^{b}}d\tau d\tau_1.
\end{equation*}

The inequality (\ref{DUA2}) with $Z_1(f,g,\phi)$ can be estimate by the same argument as the one used in the bound of $W_2(f,g,\phi)$.

Now we  proof inequality (\ref{DUA2}) with $Z_2(f,g,\phi)$. First we
make the change of variables $\tau_2=\tau-\tau_1$, $n_2=n-n_1$ to
obtain
\begin{eqnarray*}
Z_2(f,g,\phi)&=&  \sum_{n,n_2} \int_{\R^2} \dfrac{\langle n\rangle^{s}}{\langle n-n_2\rangle^{s} \langle n_2\rangle^{s}}\\
&&\times \dfrac{h(\tau-\tau_2,n-n_2)f(\tau_2,n_2)
\bar{\phi}(\tau,n)}{\langle\tau-n^2\rangle^{a}
\langle(\tau-\tau_2)-(n-n_2)^2\rangle^{b}
\langle\tau_2+n_2^2\rangle^{b}}d\tau d\tau_2.
\end{eqnarray*}
Then changing the variables $(n, \tau, n_2, \tau_2)\mapsto -(n,
\tau, n_2, \tau_2)$ we can rewrite $Z_2(f,g,\phi)$ as
\begin{eqnarray*}
Z_2(f,g,\phi)&=&  \sum_{n,n_2} \int_{\R^2} \dfrac{\langle n\rangle^{s}}{\langle n-n_2\rangle^{s} \langle n_2\rangle^{s}}\\
&&\times \dfrac{k(\tau-\tau_2,n-n_2)l(\tau_2,n_2) \bar{\psi}(\tau,n)}{\langle\tau+n^2\rangle^{a} \langle\tau-\tau_2+(n-n_2)^2\rangle^{b} \langle\tau_2-n_2^2\rangle^{b}}d\tau d\tau_2
\end{eqnarray*}
where
$$k(a,b)=h(-a,-b),\peq l(a,b)=f(-a,-b)\peq \textrm{and} \peq {\psi}(a,b)={\phi}(-a,-b).$$

Since the $L^2$-norm is preserved under the reflection operation the result follows from the estimate for $Z_1(f,g,\phi)$.
\end{enumerate}
\fim

\begin{rema}
Once the bilinear estimates in Theorem \ref{t1.1} are established,
it is a standard matter to conclude the local well-posedness
statement of Theorem \ref{t1.2}. We refer the reader to the works
\cite{KPV2}, \cite{BOP}, \cite{GTV} and \cite{LG2} for further
details.
\end{rema}

Finally we should remark that Theorem \ref{t1.4}  can be obtained
easily using Lemma \ref{l23} (see Farah \cite{LG2}). Before get to
the end of this section we state a slightly modified bilinear
estimates that will be useful in the proof of Theorem \ref{t1.5}.
\begin{coro}\label{c3.1}
Let $a,a_1,b,b_1>1/4$ and $s\geq 0$, then there exists $c>0$ depending only on $a,a_1,b,b_1,s$ such that
\begin{enumerate}
 \item [(i)]  \peq
$\left\|
uv
\right\|_{X^S_{s,-a_1}}\lesssim \left\|u\right\|_{X^S_{s,b_1}}\left\|v\right\|_{X^B_{0,b}}+ \left\|u\right\|_{X^S_{0,b_1}}\left\|v\right\|_{X^B_{s,b}}.
$
\item [(ii)]
\peq
$\left\|
u_1\bar{u}_2
\right\|_{X^B_{s,-a}}\lesssim \left\|u_1\right\|_{X^S_{s,b_1}}\left\|u_2\right\|_{X^S_{0,b_1}}+ \left\|u_1\right\|_{X^S_{0,b_1}}\left\|u_2\right\|_{X^S_{s,b_1}}.
$
\end{enumerate}
\end{coro}

\textbf{Proof. } The above estimates are direct consequence of
Theorem \ref{t1.4} and the fact that, for all $s>0$, the following
inequality holds
\begin{equation*}
 \langle\xi\rangle^{s}\leq  \langle\xi_1\rangle^{s}  +\langle\xi-\xi_1\rangle^{s}.
\end{equation*}
\fim


\section{Counterexample to the bilinear estimates}
\textbf{Proof of Theorem \ref{t1.3}}
\begin{enumerate}
\item [$(i)$] For $u\in X^S_{k,b}$ and $v\in X^B_{s,b}$ we define $f(\tau,n)\equiv \langle\tau+n^2\rangle^b\langle n\rangle^{k} \widetilde{u}(\tau,n)$ and $g(\tau,n)\equiv \langle|\tau|-\gamma(n)\rangle^b\langle n\rangle^{s} \widetilde{v}(\tau,n)$. By Lemma \ref{l3.3} the inequality (\ref{CE1}) is equivalent to
\begin{eqnarray}\label{DE}
\left\|\dfrac{\langle n\rangle^{k}} {\langle\sigma\rangle^{a}}
\sum_{n_1}\int
\dfrac{ f(\tau_1,n_1) g(\tau_2,n_2) d\tau_1}{\langle n_1\rangle^{k}\langle n_2\rangle^{s}\langle\sigma_1\rangle^b \langle\sigma_2\rangle^b}\right\|_{l^{2}_{n}L^2_{\tau}}\lesssim \|f\|_{l^{2}_{n}L^2_{\tau}}\|g\|_{l^{2}_{n}L^2_{\tau}},
\end{eqnarray}
where
\begin{equation*}
n_2=n-n_1, \peq \tau_2=\tau-\tau_1,
\end{equation*}
\centerline{$
\sigma=\tau+n^2,\peq \sigma_1=\tau_1+n_1^2,\peq \sigma_2=|\tau_2|-n_2^2.
$}

For $N\in\Z$ define

\centerline{$f_N(\tau,n)=a_n\chi((\tau+n^2)/2)$, with $a_n=\left\{
\begin{array}{l}
1, \peq n=0, \\
0, \peq \textrm{elsewhere.}
\end{array} \right.$}
and

\centerline{$g_N(\tau,n)=b_n\chi((\tau+n^2)/2)$, with $b_n=\left\{
\begin{array}{l}
1, \peq n=N, \\
0, \peq \textrm{elsewhere.}
\end{array} \right.$}
where $\chi(\cdot)$ denotes the characteristic function of the interval $[-1,1]$. Thus
$$a_{n_1}b_{n-n_1}\neq 0 \peq \textrm{if and only if} \peq n_1=0 \peq \textrm{and} \peq n=N$$
and consequently for $N$ large
\begin{eqnarray*}
\int\!\!\chi((\tau_1+n_1^2)/2)\chi((\tau-\tau_1+(n-n_1)^2)/2) &\!\!\!\!\gtrsim\!\!\!\!& \chi((\tau+(n-n_1)^2+n_1^2))\\
&\!\!\!\!\gtrsim\!\!\!\!& \chi((\tau+N^2)).
\end{eqnarray*}
Therefore, using the fact that $||\tau_2|-n_2^2|\leq |\tau_2+n_2^2|$, inequality (\ref{DE}) implies
\begin{eqnarray*}
1 \peq \gtrsim \peq \|N^{k-s}\chi((\tau+N^2))\|_{L^2_{\tau}}\peq \gtrsim \peq N^{k-s}.
\end{eqnarray*}

Letting $N\rightarrow\infty$, this inequality is possible only when $k\leq s$.
\item [$(ii)$] Now we define

\centerline{$f_N(\tau,n)=a_n\chi((\tau+n^2)/2)$, with $a_n=\left\{
\begin{array}{l}
1, \peq n=-N, \\
0, \peq \textrm{elsewhere.}
\end{array} \right.$}
and

\centerline{$g_N(\tau,n)=b_n\chi((\tau-n^2)/2)$, with $b_n=\left\{
\begin{array}{l}
1, \peq n=N, \\
0, \peq \textrm{elsewhere.}
\end{array} \right.$}
Then
$$a_{n_1}b_{n-n_1}\neq 0 \peq \textrm{if and only if} \peq n_1=0 \peq \textrm{and} \peq n=N$$
and for $N$ large
\begin{eqnarray*}
\int\!\!\chi((\tau_1+n_1^2)/2)\chi((\tau-\tau_1-(n-n_1)^2)/2) &\!\!\!\!\gtrsim\!\!\!\!& \chi((\tau+n^2-2nn_1))\\
&\!\!\!\!\gtrsim\!\!\!\!& \chi((\tau)).
\end{eqnarray*}
Therefore, using the fact that $||\tau_2|-n_2^2|\leq |\tau_2-n_2^2|$, inequality (\ref{DE}) implies
\begin{eqnarray*}
1 \peq \gtrsim \peq \|N^{-(k+s)}\chi((\tau))\|_{L^2_{\tau}}\peq \gtrsim \peq N^{-(k+s)}.
\end{eqnarray*}

Letting $N\rightarrow\infty$, this inequality is possible only when $k+s\geq 0$.
\item [$(iii)$] For $u_1\in X^S_{k,b}$ and $u_2\in X^S_{k,b}$ we define $f(\tau,n)\equiv \langle\tau+n^2\rangle^b\langle n\rangle^{k} \widetilde{u}_1(\tau,\xi)$ and $g(\tau,n)\equiv \langle\tau+n^2\rangle^b\langle n\rangle^{k} \widetilde{u}_2(\tau,\xi)$. By Lemma \ref{l3.3} the inequality (\ref{CE2}) is equivalent to
\begin{eqnarray}\label{DE2}
\left\|\dfrac{\langle n\rangle^{s}} {\langle\sigma\rangle^{a}}
\sum_{n_1}\int
\dfrac{ f(\tau_1,n_1) h(\tau_2,n_2) d\tau_1}{\langle n_1\rangle^{k}\langle n_2\rangle^{k}\langle\sigma_1\rangle^b \langle\sigma_2\rangle^b}\right\|_{l^{2}_{n}L^2_{\tau}}\lesssim \|f\|_{l^{2}_{n}L^2_{\tau}}\|g\|_{l^{2}_{n}L^2_{\tau}},
\end{eqnarray}
where
\begin{equation*}
h(\tau_2,n_2)=\bar{g}(-\tau_2,-n_2), \peq n_2=n-n_1, \peq \tau_2=\tau-\tau_1,
\end{equation*}
\centerline{$
\sigma=|\tau|-n^2,\peq \sigma_1=\tau_1+n_1^2,\peq \sigma_2=\tau_2-n_2^2.
$}

For $N\in\Z$ define

\centerline{$f_N(\tau,n)=a_n\chi((\tau+n^2)/2)$, with $a_n=\left\{
\begin{array}{l}
1, \peq n=N, \\
0, \peq \textrm{elsewhere.}
\end{array} \right.$}
and

\centerline{$h_N(\tau,n)=b_n\chi((\tau-n^2)/2)$, with $b_n=\left\{
\begin{array}{l}
1, \peq n=0, \\
0, \peq \textrm{elsewhere.}
\end{array} \right.$}
where $\chi(\cdot)$ denotes the characteristic function of the interval $[-1,1]$. Thus
$$a_{n_1}b_{n-n_1}\neq 0 \peq \textrm{if and only if} \peq n_1=N \peq \textrm{and} \peq n=N$$
and
\begin{eqnarray*}
\int\!\!\chi((\tau_1+n_1^2)/2)\chi((\tau-\tau_1-(n-n_1)^2)/2) &\!\!\!\!\gtrsim\!\!\!\!& \chi((\tau-(n-n_1)^2+n_1^2))\\
&\!\!\!\!\gtrsim\!\!\!\!& \chi((\tau+N^2)).
\end{eqnarray*}
Therefore, using the fact that $||\tau|-n^2|\leq |\tau+n^2|$, inequality (\ref{DE2}) implies
\begin{eqnarray*}
1 \peq \gtrsim \peq \|N^{s-k}\chi((\tau+N^2))\|_{L^2_{\tau}}\peq \gtrsim \peq N^{s-k}.
\end{eqnarray*}

Letting $N\rightarrow\infty$, this inequality is possible only when $s\leq k$.\\
\fim
\end{enumerate}

\section{Global Well-posedness}

We divide our analysis in two cases. The proof of Theorem \ref{t1.5}
for $s=0$ follows the same  lines as in Farah \cite{LG2} Theorem
$1.4$. For the convenience of the reader we repeat the proof of this
case below. The case $s>0$ can be
proved using the arguments introduced by Bourgain \cite{B5} for the
Schr\"odinger equation and further developed by Angulo \textit{et
al.} \cite{AMP1} for the Schr\"odinger-Benjamin-Ono system.\\

\textbf{Proof of Theorem \ref{t1.5}.}\\

 \textit{\underline{Case $s=0$}}:\\

 Let $({u_0},{v_0},v_1)\in L_{per}^2([0,1])\times L_{per}^2([0,1])\times H_{per}^{-1}([0,1])$ and $0<T\leq 1$. Based on the integral formulation (\ref{INT4}), we define the integral operators
\begin{eqnarray}\label{INT5}
\begin{small}
\begin{split}
G_T^S(u,v)(t)=& \theta_T(t)U(t)u_0-i\theta_T(t)\int_{0}^{t}U(t-t')(\alpha vu)(t')dt'\\
G_T^B(u,v)(t)=& \theta_T(t)\left(V_c(t)v_0+V_s(t)(v_1)_x\right)+\theta_T(t) \int_{0}^{t}V_s(t-t')(\beta|u|^2)_{xx}(t')dt'.
\end{split}
\end{small}
\end{eqnarray} Therefore, applying Lemmas \ref{l21}-\ref{l22} and Theorem \ref{t1.3}, we obtain
\begin{eqnarray}\label{C2}
\begin{split}
\|G_T^S(u,v)\|_{X^S_{0,b_1}}&\leq cT^{1/2-b_1}\|u_0\|_{L_{per}^2}+cT^{1-(a_1+b_1)}\left\|uv\right\|_{X^S_{0,-a_1}}\\
&\leq cT^{1/2-b_1}\|u_0\|_{L_{per}^2}+cT^{1-(a_1+b_1)}\left\|u\right\|_{X^S_{0,b_1}} \left\|v\right\|_{X^B_{0,b}},\\
\|G_T^B(u,v)\|_{X^B_{0,b}}&\leq cT^{1/2-b}\|v_0,v_1\|_{\mathfrak{B}}+cT^{1-(a+b)}\left\|u\bar{u}\right\|_{X^B_{0,-a}}\\
&\leq cT^{1/2-b}\|v_0,v_1\|_{\mathfrak{B}}+cT^{1-(a+b)}\left\|u\right\|_{X^S_{0,b_1}}^2
\end{split}
\end{eqnarray}
and also
\begin{eqnarray}\label{C3}
\begin{split}
\|G_T^S(u,v)-G_T^S(z,w)\|_{X^S_{0,b_1}}&\leq cT^{1-(a_1+b_1)}\left(\left\|u\right\|_{X^S_{0,b_1}} \left\|v-w\right\|_{X^B_{0,b}}\right.\\
&\left.+\left\|u-z\right\|_{X^S_{0,b_1}} \left\|w\right\|_{X^B_{0,b}}\right),\\
\|G_T^B(u,v)-G_T^B(z,w)\|_{X^B_{0,b}}&\leq cT^{1-(a+b)}\left(\left\|u\right\|_{X^S_{0,b_1}} + \left\|z\right\|_{X^S_{0,b_1}}\right)\\
&\times \left\|u-z\right\|_{X^S_{0,b_1}}.\\
\end{split}
\end{eqnarray}

We define
\begin{eqnarray*}
X_{0,b_1}^S(d_1)&=&\left\{u\in X_{0,b_1}^S:\|u\|_{X_{0,b_1}^S}\leq d_1\right\},\\
X_{0,b}^B(d)&=&\left\{v\in X_{0,b}^B:\|v\|_{X_{0,b}^B}\leq d\right\},
\end{eqnarray*}
where $d_1=2cT^{1/2-b_1}\|u_0\|_{L_{per}^2}$ and $d=2cT^{1/2-b}\|v_0,v_1\|_{\mathfrak{B}}$.\\

For $(G_T^S, G_T^B)$ to be a contraction in $X_{0,b_1}^S(d_1)\times X_{0,b}^B(d)$ it needs to satisfy
\begin{eqnarray}\label{GC1}
d_1/2+cT^{1-(a_1+b_1)}d_1d\leq d_1
\Leftrightarrow T^{3/2-(a_1+b_1+b)}\|v_0,v_1\|_{\mathfrak{B}}\lesssim 1,
\end{eqnarray}
\begin{eqnarray}\label{GC2}
d/2+cT^{1-(a+b)}d_1^2\leq d
\Leftrightarrow T^{3/2-(a+2b_1)}\|u_0\|_{L_{per}^2}^2\lesssim \|v_0,v_1\|_{\mathfrak{B}},
\end{eqnarray}
\begin{eqnarray}\label{GC3}
2cT^{1-(a+b)}d_1\leq 1/2
\Leftrightarrow T^{3/2-(a+b+b_1)}\|u_0\|_{L_{per}^2}\lesssim 1,
\end{eqnarray}
\begin{eqnarray}\label{GC5}
2cT^{1-(a_1+b_1)}d_1\leq 1/2
\Leftrightarrow T^{3/2-(a_1+2b_1)}\|u_0\|_{L_{per}^2}\lesssim 1.
\end{eqnarray}

Therefore, we conclude that there exists a solution $(u,v)\in X_{0,b_1}^S\times X_{0,b}^B$ satisfying
\begin{eqnarray}\label{GC}
\|u\|_{X_{0,b_1}^{S,[0,T]}}\leq 2cT^{1/2-b_1}\|u_0\|_{L_{per}^2} \textrm{\peq and \peq}
\|v\|_{X_{0,b}^{B,[0,T]}}\leq 2cT^{1/2-b}\|v_0,v_1\|_{\mathfrak{B}}.
\end{eqnarray}

On the other hand, applying Lemmas \ref{l21}-\ref{l22} we have that, in fact, $(u,v)\in C([0,T]:L^2)\times C([0,T]:L^2)$. Moreover, since the $L^2$-norm of $u$ is conserved by the flow we have $\|u(T)\|_{L_{per}^2}=\|u_0\|_{L_{per}^2}$.

Now, we need to control the growth of $\|v(t)\|_{\mathfrak{B}}$ in each time step. If, for all $t>0$, $\|v(t)\|_{\mathfrak{B}}\lesssim \|u_0\|_{L_{per}^2}^2$ we can repeat the local well-posedness argument and extend the solution globally in time. Thus, without loss of generality, we suppose that after some number of iterations we reach a time $t^{\ast}>0$ where $\|v(t^{\ast})\|_{\mathfrak{B}}\gg \|u_0\|_{L_{per}^2}^2$.

Hence, since $0<T\leq 1$, condition (\ref{GC2}) is automatically satisfied and conditions (\ref{GC1})-(\ref{GC5}) imply that we can select a time increment of size
\begin{eqnarray}\label{GC4}
T\sim \|v(t^{\ast})\|_{\mathfrak{B}}^{-1/(3/2-(a_1+b_1+b))}.
\end{eqnarray}

Therefore, applying Lemmas \ref{l21}($b$)-\ref{l22}($b$) to $v=G_T^B(u,v)$ we have
\begin{eqnarray*}
\|v(t^{\ast}+T)\|_{\mathfrak{B}}\leq \|v(t^{\ast})\|_{\mathfrak{B}}+cT^{3/2-(a+2b_1)}(\|u_0\|_{L_{per}^2}^2+1).
\end{eqnarray*}

Thus, we can carry out $m$ iterations on time intervals, each of length (\ref{GC4}), before the quantity $\|v(t)\|_{\mathfrak{B}}$ doubles, where $m$ is given by
\begin{eqnarray*}
mT^{3/2-(a+2b_1)}(\|u_0\|_{L_{per}^2}^2+1)\sim \|v(t^{\ast})\|_{\mathfrak{B}}.
\end{eqnarray*}

The total time of existence we obtain after these $m$ iterations is
\begin{eqnarray*}
\Delta T=mT&\sim& \dfrac{\|v(t^{\ast})\|_{\mathfrak{B}}}{T^{1/2-(a+2b_1)} (\|u_0\|_{L_{per}^2}^2+1)}\\
&\sim& \dfrac{\|v(t^{\ast})\|_{\mathfrak{B}}} {\|v(t^{\ast})\|_{\mathfrak{B}}^{-(1/2-(a+2b_1))/(3/2-(a_1+b_1+b))} (\|u_0\|_{L_{per}^2}^2+1)}.
\end{eqnarray*}

Taking $a,b,a_1,b_1$ such that
\begin{eqnarray*}
 \dfrac{a+2b_1-1/2}{(3/2-(a_1+b_1+b))}=1,
\end{eqnarray*}
(for instance, $a=b=a_1=b_1=1/3$) we have that $\Delta T$  depends
only on $\|u_0\|_{L_{per}^2}$, which is conserved by the flow. Hence
we can repeat this entire argument and extend the solution $(u,v)$
globally in time.

Moreover, since in each step of time $\Delta T$ the size  of
$\|v(t)\|_{\mathfrak{B}}$ will at most double it is easy to see
that, for all $\widetilde{T}>0$
\begin{eqnarray}\label{GGV5}
 \|v(\widetilde{T})\|_{\mathfrak{B}}\lesssim \exp{((\ln{2})\|u_0\|_{L_{per}^2}^2\widetilde{T})} \max{\{\|v_0,v_1\|_{\mathfrak{B}},\|u_0\|_{L_{per}^2}\}}.
\end{eqnarray}

\textit{\underline{Case $s>0$}}:\\

Let $(u_0,v_0,v_1) \in H_{per}^s\times H_{per}^s\times
H_{per}^{s-1}$.  By the previous case, there exists a global
solution $(u,v)\in C([0,+\infty);L_{per}^2)\times
C([0,+\infty);L_{per}^2)$. Moreover, $(u,v)$ is a solution of the
integral equation (\ref{INT5}) in the time interval $[0,\Delta T]$,
with $\Delta T \sim \dfrac{1}{\|u_0\|_{L_{per}^2}^2+1}$, satisfying
\begin{equation}\label{s0}
\max\left\{\|u\|_{X_{0,1/3}^{S,[0,\Delta T]}},\|v\|_{X_{0,1/3}^{B,[0,\Delta T]}}\right\}\lesssim C(\|u_0\|_{L_{per}^2},\|v_0,v_1\|_{\mathfrak{B}}),
\end{equation}
where the constant
$C(\|u_0\|_{L_{per}^2},\|v_0,v_1\|_{\mathfrak{B}})>0$  depends only
on $\|u_0\|_{L_{per}^2}$ and $\|v_0,v_1\|_{\mathfrak{B}}$.

We claim that the solution $(u,v)$, in fact, belongs to $X_{s,1/3}^{S,[0,T_0]}\times X_{s,1/3}^{B,[0,T_0]}$ for all $0 < T_0 \leq \Delta T$. Indeed, applying Lemmas \ref{l21}-\ref{l22} and Corollary \ref{c3.1} with $a=b=a_1=b_1=1/3$, we obtain
\begin{equation}\label{s01}
\|u\|_{X^{S,[0,T_0]}_{s,1/3}}\lesssim \|u_0\|_{H_{per}^s}+T_0^{1/3}\left(\left\|u\right\|_{X^{S,[0,T_0]}_{s,1/3}} \left\|v\right\|_{X^{B,[0,T_0]}_{0,1/3}}+ \left\|u\right\|_{X^{S,[0,T_0]}_{0,1/3}}\left\|v\right\|_{X^{B,[0,T_0]}_{s,1/3}}\right)
\end{equation}
and
\begin{equation}\label{s02}
\|v\|_{X^{B,[0,T_0]}_{s,1/3}}\lesssim \|v_0,v_1\|_{\mathfrak{B}^s}+ T_0^{1/3}\left(\left\|u\right\|_{X^{B,[0,T_0]}_{s,1/3}} \left\|u\right\|_{X^{B,[0,T_0]}_{0,1/3}}\right),
\end{equation}
where $0 < T_0 \leq \Delta T$. Inserting the inequality (\ref{s02}) into (\ref{s01}) and using (\ref{s0}) we conclude
\begin{eqnarray*}
\|u\|_{X^{S,[0,T_0]}_{s,1/3}}&\lesssim& \|u_0\|_{H_{per}^s}+C(\|u_0\|_{L_{per}^2},\|v_0,v_1\|_{\mathfrak{B}}) \|v_0,v_1\|_{\mathfrak{B}^s}\\ &&+T_0^{1/3}C(\|u_0\|_{L_{per}^2},\|v_0,v_1\|_{\mathfrak{B}}) \left\|u\right\|_{X^{S,[0,T_0]}_{s,1/3}}.
\end{eqnarray*}

Set
\begin{eqnarray*}
T_0 \sim \dfrac{1} {\left(1+C(\|u_0\|_{L_{per}^2},\|v_0,v_1\|_{\mathfrak{B}})\right)^3}.
\end{eqnarray*}

Hence, from the choice of $T_0$, we deduce the following \textit{a
priori} estimates
\begin{equation*}
\|u\|_{X^{S,[0,T_0]}_{s,1/3}}\lesssim \|u_0\|_{H_{per}^s}+C(\|u_0\|_{L_{per}^2},\|v_0,v_1\|_{\mathfrak{B}}) \|v_0,v_1\|_{\mathfrak{B}^s}
\end{equation*}
and
\begin{equation*}
\|v\|_{X^{B,[0,T_0]}_{s,1/3}}\lesssim \|v_0,v_1\|_{\mathfrak{B}^s}+C(\|u_0\|_{L_{per}^2},\|v_0,v_1\|_{\mathfrak{B}}) \left(\|v_0,v_1\|_{\mathfrak{B}^s}+\|u_0\|_{L_{per}^2}\right).
\end{equation*}

Thus, applying Lemmas \ref{l21}-\ref{l22} we get that $(u,v)\in
C([0,T_0];H_{per}^s)\times C([0,T_0];H_{per}^s)$.  The preceding
statement remains valid for any bounded interval $[0,T]$, since
$T_0$ depends only on $\|u_0\|_{L_{per}^2}$ and
$\|v_0,v_1\|_{\mathfrak{B}}$ and we can iterate the above argument a
finite number of times to deduce that
\begin{equation*}
(u,v)\in C([0,T];H_{per}^s)\times C([0,T];H_{per}^s), \textrm{ for all $T>0$},
\end{equation*}
which completes the proof of Theorem \ref{t1.5}.\\
\fim
\section{Stability of periodic traveling waves}

As we said in the Introduction, here we will consider system
\eqref{SB} with $\alpha=\beta=-1$, that is, we consider the system
\begin{eqnarray}\label{s1}
\left\{
\begin{array}{l}
iu_{t}+u_{xx}+ uv=0, \\
v_{tt}-v_{xx}+v_{xxxx}+(|u|^2)_{xx}=0,\\
\end{array} \right.
\end{eqnarray}
and look for traveling waves of the form
\begin{equation}\label{s2}
u(x,t)=e^{i\omega t} \psi_{\omega}(x), \qquad v(x,t)=\phi_\omega(x),
\end{equation}
where $\omega$ is a real parameter (to be determined later) and
$\psi_\omega, \phi_\omega:\mathbb{R}\rightarrow\mathbb{R}$ are
smooth periodic functions with the same fixed period $L>0$. Then,
substituting \eqref{s2} into \eqref{s1}; integrating twice the
second equation in the obtained system and assuming that the
integration constants are zero, we obtain the system
\begin{eqnarray}\label{s4}
\left\{
\begin{array}{l}
\psi_\omega''-\omega\psi_\omega+\psi_\omega\phi_\omega=0, \\
\phi_\omega''-\phi_\omega+\psi_\omega^2=0.\\
\end{array} \right.
\end{eqnarray}
In order to solve system \eqref{s4} we assume $\omega=1$ and
$\psi_1=\phi_1$, so that system \eqref{s4} reduces to a single
ordinary differential equation, namely,
\begin{equation}\label{s5}
\psi_1''-\psi_1+\psi_1^2=0.
\end{equation}

As we will see later in our stability analysis, it is necessary to
construct a smooth branch of periodic wave solutions (depending on
$\omega$) passing through solution $\psi_1$ of \eqref{s5}. Then, we
will consider the family of equations
\begin{equation}\label{s6}
\psi_\omega''-\omega\psi_\omega+\psi_\omega^2=0,
\end{equation}
so that at $\omega=1$ we obtain a solution for \eqref{s5}.

\subsection{Existence of traveling waves}

Along this subsection, we review the theory of finding solutions for \eqref{s6}.
Indeed,  equation \eqref{s6} can be solved by using the standard
{\it direct integration method} (for details we refer to \cite{AL}).
As a matter of fact, equation \eqref{s6} has a \emph{strictly
positive} solution of the form
\begin{equation}   \label{s7}
\displaystyle\psi_{\omega}(x)=
\beta_2+(\beta_3-\beta_2)cn^2\left(\sqrt{\frac{\beta_3-\beta_1}{6}}x;k\right),
\qquad k^2=\frac{\beta_3-\beta_2}{\beta_3-\beta_1},
\end{equation}
where $cn(\cdot;k)$ denotes the Jacobian elliptic function of {\it
cnoidal} type, $k$ is the elliptic modulus and
$\beta_1,\beta_2,\beta_3$ are real constants satisfying
\begin{equation}\label{s7.1}
\frac{3\omega}{2}=\sum_{i=1}^3\beta_i, \qquad
0=\sum_{i<j}\beta_i\beta_j, \qquad \beta_1\beta_2\beta_3=3A_\psi,
\end{equation}
where $A_\psi$ is an integration constant. Moreover, it must be the
case that
$$
\beta_1<0<\beta_2<\omega<\beta_3<\frac{3\omega}{2}.
$$

The first question concerning solution \eqref{s7} is the following:
Fixed $L>0$, can we choose $\beta_1,\beta_2,\beta_3$ such that
solution \eqref{s7} has fundamental period $L$? The answer is yes.
To prove so, one first note since $cn^2(\cdot;k)$ has fundamental
period $2K(k)$, where $K$ is the complete elliptic integral of the
first kind defined by (see e.g., \cite{BF})
$$
K(k)=\int_0^1\;\frac{dt}{\sqrt{(1-t^2)(1-k^2t^2)}},
$$
function $\psi_\omega$ given in \eqref{s7} has fundamental period
\begin{equation}\label{s8}
\displaystyle
T_{\psi_{\omega}}=\frac{2\sqrt{6}}{\sqrt{\beta_3-\beta_1}}K(k).
\end{equation}

Next we observe that $T_{\psi_{\omega}}$ can be rewritten as a
function depending only on $\beta_2$ (and $\omega>0$ fixed). In
fact, by defining $\omega_0=\omega/2$, we readily see from
\eqref{s7.1} that
\begin{equation}\label{s9}
\displaystyle
T_{\psi_{\omega}}(\beta_2;\omega_0)=\frac{2\sqrt{6}}{\sqrt{\rho(\beta_2;\omega_0)}}K(k(\beta_2;\omega_0)),
\end{equation}
where
\begin{equation}\label{s10}
\rho(\beta_2;\omega_0)=\sqrt{9\omega_0^2-3\beta_2^2+6\omega_0\beta_2},
\qquad
k^2(\beta_2;\omega_0)=\frac{1}{2}+\frac{3(\omega_0-\beta_2)}{2\rho(\beta_2;\omega_0)}.
\end{equation}

Moreover, from \eqref{s9} it is easy to see that
$T_{\psi_{\omega}}\rightarrow +\infty$, as $\beta_2\rightarrow0$ and
$T_{\psi_{\omega}}\rightarrow \sqrt2\pi/\sqrt{\omega_0}$, as
$\beta_2\rightarrow 2\omega_0$. Since the function
$\beta_2\in(0,2\omega_0)\rightarrow
T_{\psi_{\omega}}(\beta_2;\omega_0)$ is strictly decreasing (this
will be proved in the next theorem) we see that, fixed $L>0$ and
choosing $\omega_0>2\pi^2/L^2$, there exists a unique $\beta_2\equiv
\beta_2(\omega_0)\in(0,2\omega_0)$ such that the corresponding
cnoidal wave given by \eqref{s7} has fundamental period
$T_{\psi_{\omega}}(\beta_2;\omega_0)=L$.

In supplement to the above analysis, fixed $L>0$, we can construct a
smooth curve (depending on $\omega$) of cnoidal waves solutions for
\eqref{s6} such that each one of its elements have fundamental
period $L$. This is the content of the next theorem.

\begin{theo}  \label{exist}
Let $L>2\pi$ be fixed. Choose arbitrarily $\omega_0>2\pi^2/L^2$ and
consider the unique $\beta_{2,0}=\beta_2(\omega_0)\in (0,2\omega_0)$
such that
$$
L=\frac{2\sqrt6}{\sqrt{\rho(\beta_{2,0};\omega_0)}}K(k(\beta_{2,0};\omega_0)).
$$
Then,
\begin{itemize}
  \item[(i)] there exist an interval $J_1(\omega_0)$ around
  $\omega_0$, an interval $J_2(\beta_{2,0})$ around $\beta_{2,0}$
  and a unique smooth function $\Lambda:  J_1(\omega_0)\rightarrow
  J_2(\beta_{2,0})$ such that $\Lambda(\omega_0)=\beta_{2,0}$ and
  $$
  L=
  \frac{2\sqrt6}{\sqrt{\rho(\beta_{2};\eta)}}K(k(\beta_{2};\eta)),
  $$
  where $\eta\in J_1(\omega_0), \beta_2=\Lambda(\eta)$ and
  $k(\beta_2;\eta), \rho(\beta_2;\eta)$ are defined in \eqref{s10}
  with $\omega_0$ replaced with $\eta$. Moreover, the interval
  $J_1(\omega_0)$ can be chosen to be the interval $\mathcal{I}=(2\pi^2/L^2,
  +\infty)$ and the modulus $k=k(\eta)$, where
  \begin{equation}\label{sk}
    k^2(\eta):=\frac{1}{2}+\frac{3(\eta-\Lambda(\eta))}{2\rho(\Lambda(\eta);\eta)},
  \end{equation}
   is a strictly increasing
  function (on the parameter  $\eta$).

  \item[(ii)] For $\omega\in (4\pi^2/L^2, +\infty)$ and $\eta(\omega)=\omega/2$, the cnoidal
  wave solution $\psi_\omega(\cdot)=\psi_{\eta(\omega)}(\cdot; \beta_2(\eta(\omega)))$ has fundamental
  period $L$ and satisfies \eqref{s6}. In addition, the mapping
  $$
  \omega\in\left(\frac{4\pi^2}{L^2},+\infty\right) \mapsto \psi_\omega\in
  H^k_{per}([0,L]), \quad k=0,1,\ldots
  $$
  is a smooth function.
\end{itemize}
\end{theo}
\textbf{Sketch of the proof.} The proof is an application of the
Implicit Function Theorem. Here we give only the main steps (for
details see \cite{AL}). Define $\Omega=\{(\beta_2,\eta)\in
\mathbb{R}^2; \,\, \eta>2\pi^2/L^2, \beta_2\in (0,2\eta)\}$ and
$\Gamma: \Omega\rightarrow\mathbb{R}$ by
$$
\Gamma(\beta_2,\eta)=
\frac{2\sqrt6}{\sqrt{\rho(\beta_{2};\eta)}}K(k(\beta_{2};\eta))-L.
$$

By our assumptions, we have $\Gamma(\beta_{2,0},\omega_0)=0$.
Moreover, taking into account the properties of the complete
elliptic integrals and the definitions of $k$ and $\rho$ one infers
that $\partial\Gamma/\partial\beta_2<0$ for all $(\beta_2,\eta)\in
\Omega$. So, an application of the Implicit Function Theorem gives
us the desired. The fact that $J_1(\omega_0)$ can be chosen to be
$\mathcal{I}$ follows from the fact that $\omega_0$ can be
arbitrarily chosen in $\mathcal{I}$ and the uniqueness of the
function arising in the Implicit Function Theorem.

To see that $k(\eta)$ is a strictly increasing function one just
take the derivative with respect to $\eta$ in \eqref{sk} and note
that $dk/d\eta>0$. \\
\fim

\begin{rema}
We have assumed $L>2\pi$ in Theorem \ref{exist} because we want to
get a smooth curve of cnoidal  wave (defined in an open interval)
passing through $\omega=1$. Otherwise, that is, if $L\leq2\pi$ then
such a curve does not exist.
\end{rema}

\subsection{Spectral Analysis}

To obtain our stability results, we will use the Grillakis, Shatah
and Strauss theory \cite{GSS}. As it is well-known in such approach
we need to study the spectrum of some linearized operators.

First, we note that introducing a new variable $w$ defined
by $v_t=w_x$, system \eqref{s1} can be written as an Hamiltonian
system of the form
\begin{equation}\label{sa1}
\frac{d}{dt}U(t)=J\mathcal{E}'(U(t)),
\end{equation}
where $U=(P,v,Q,w)$, $P=\mbox{Re}(u)$, $Q=\mbox{Im}(u)$, $J$ is the
skew-symmetric matrix
\begin{equation}   \label{sa2}
J= \left( \begin{array}{cccc}
0 & 0 & 1/2 & 0    \\
0 & 0 & 0 & \partial_x   \\
-1/2 & 0 & 0 & 0    \\
0 & \partial_x & 0 & 0
\end{array}
\right)
\end{equation}
and $\mathcal{E}$ is the energy functional given by
\begin{equation}\label{sa3}
\mathcal{E}(U)=\int_0^L
\left\{P_x^2+Q_x^2+\frac{v_x^2}{2}+\frac{v^2}{2}+\frac{w^2}{2}-v(P^2+Q^2)\right\}dx.
\end{equation}

Next we will consider the linearized operator we need to study. We
first remind that system \eqref{s1} preserves the $L^2$ norm of $u$
and  so, in the above notation,
$$
\mathcal{F}(U)=\int_0^L \{P^2+Q^2\}dx
$$
is a conserved quantity of system \eqref{s1}.

To simplify our exposition, we  denote
$\Psi_\omega=(\psi_\omega,\psi_\omega,0,0)$, where $\psi_\omega$ is
a cnoidal wave given in Theorem \ref{exist}. By direct computation
we see that $\Psi_\omega$ is a critical point of the functional
$\mathcal{E}+\omega\mathcal{F}$ at $\omega=1$, that is,
\begin{equation}\label{critical}
\mathcal{E}'(\Psi_1)+\mathcal{F}'(\Psi_1)=0.
\end{equation}

Now consider the operator
\begin{equation}\label{sa4}
\mathcal{A}:=\mathcal{E}''(\Psi_1)+\mathcal{F}''(\Psi_1)=
\left(\begin{array}{cccc}
\mathcal{A}_{R} & 0\\\\
0 & \mathcal{A}_I
\end{array}\right),
\end{equation}
where $\mathcal{A}_R$ and $\mathcal{A}_I$ are the self-adjoint
$2\times 2$ matrix differential operators defined by
\begin{equation} \label{sa5}
\displaystyle\mathcal{A}_{R}=\left(
\begin{array}{cccc}
\displaystyle 2(-\partial_x^2+1-\psi_1) & -2\psi_1\\\\
-2\psi_1 & \displaystyle -\partial_x^2+1
\end{array}\right)
\end{equation}
 and
\begin{equation}  \label{sa6}
\displaystyle\mathcal{A}_{I}=\left(
\begin{array}{cccc}
\displaystyle2(-\partial_x^2+1-\psi_1) & 0\\\\
0 & 1
\end{array}\right).
\end{equation}

Let us study the spectrum of operator $\mathcal{A}$. In what
follows, we use the notation $\sigma(\mathcal{L})$ to represent the
spectrum of the linear operator $\mathcal{L}$. We first remind that
if $\sigma_{ess}(\mathcal{L})$ and $\sigma_{disc}(\mathcal{L})$
denote, respectively, the essential and discrete spectra of
$\mathcal{L}$, then
$\sigma(\mathcal{L})=\sigma_{ess}(\mathcal{L})\cup\sigma_{disc}(\mathcal{L})$.

To begin our analysis, we observe that since $\mathcal{A}$ is a
diagonal operator we have
$\sigma(\mathcal{A})=\sigma(\mathcal{A}_R)\cup\sigma(\mathcal{A}_I)$.
Moreover, sice $\mathcal{A}$ has a compact resolvent we obtain
$\sigma(\mathcal{A})=\sigma_{disc}(\mathcal{A})$ (see e.g.,
\cite{RSIV})

Before studying the spectrum of operators $\mathcal{A}_R$ and
$\mathcal{A}_I$, we recall the following lemma

\begin{lemm} \label{lemmaAng}
Let $\psi=\psi_{1}$ be the cnoidal wave given by Theorem
\ref{exist}. Then the following spectral properties hold
\begin{itemize}
    \item [(i)] Operator
    $$
    \mathcal{L}_{1}:=-\partial_x^2+1-2\psi
    $$
    defined in $L^2_{per}([0,L])$ with domain $H_{per}^2([0,L])$ has exactly one negative
    eigenvalue which is simple;  zero is an eigenvalue which is
    simple with eigenfunction $\psi'$. Moreover, the remainder
    of the spectrum is constituted by a discrete set of
    eigenvalues.
    \item [(ii)] Operator
     $$
    \mathcal{L}_{2}:=-\partial_x^2+1-\psi
    $$
    defined in $L^2_{per}([0,L])$ with domain $H_{per}^2([0,L])$
    has no negative eigenvalues; zero is an eigenvalue, simple
    with eigenfunction $\psi$. Moreover, the remainder
    of the spectrum is constituted by a discrete set of
    eigenvalues.
\end{itemize}
\end{lemm}
\textbf{Proof.} For the first part, see Theorem 4.1 in \cite{AL}.
The second part follows immediately from Floquet's theory. Indeed,
in view of \eqref{s5} we have that $0$ is an eigenvalue for
$\mathcal{L}_{2}$ with eigenfunction $\psi$. Moreover, since $\psi$
has no zeros in the interval $[0,L]$, 0 must be the first
eigenvalue (see e.g. \cite[Chapter 3]{EASTHAM}).\\
\fim

With Lemma \ref{lemmaAng} at hands, we are able to prove some
spectral properties for operators $\mathcal{A}_R$ and
$\mathcal{A}_I$.

\begin{theo}   \label{teoeigenR}
Let $\psi=\psi_{1}$ be the \textit{cnoidal} wave solution given by
Theorem \ref{exist}. Then,
\begin{itemize}
\item[(i)] operator $\mathcal{A}_{R}$ in \eqref{sa5}
defined in $L_{per}^2([0,L])\times L_{per}^2([0,L])$ with domain
$H_{per}^2([0,L])\times H_{per}^2([0,L])$ has  its first three
eigenvalues simple, being the eigenvalue zero the second one with
eigenfunction $(\psi',  \psi')$. Moreover, the remainder of the
spectrum is constituted by a discrete set of eigenvalues.

  \item [(ii)] Operator $\mathcal{A}_{I}$ in \eqref{sa6}
  defined in $L_{per}^2([0,L])\times L_{per}^2([0,L])$ with domain
 $H_{per}^2([0,L])\times L_{per}^2([0,L])$ has no negative
eigenvalues; zero is the first eigenvalue which is simple with
eigenfunction $(\psi,0)$. Moreover, the remainder of the spectrum is
constituted by a discrete set of eigen\-values.
\end{itemize}
\end{theo}
\textbf{Proof.} (i) First we observe that from \eqref{s6} it is easy
to see that zero is an eigenvalue with eigenfunction $(\psi',
\psi')$. Now we consider the quadratic form associated with
$\mathcal{A}_R$. Let  $Y=H^1_{per}([0,L])\times H^1_{per}([0,L])$.
Then, for $(f,g)\in Y$,
\begin{equation} \label{sa7}
\begin{array}{lll}
Q_R(f,g)&:=&\displaystyle \langle\mathcal{A}_R(f,g),(f,g)\rangle\\\\
&=&\displaystyle\int_{0}^{L}\{ 2(-\partial_x^2+1-\psi)f^2-4\psi fg+(-\partial_x^2+1)g^2 \}\,dx\\\\
&=& 2\langle\mathcal{L}_1f,f\rangle + \langle\mathcal{L}_1g,g\rangle
+ \displaystyle 2\int_0^L\psi(f-g)^2\, dx.
\end{array}
\end{equation}

In order to prove that $\mathcal{A}_R$ has at least one negative eigenvalue,
let us prove that there exists a pair $(f,g)\in Y$ such that
$Q_R(f,g)<0$. Indeed, from Lemma \ref{lemmaAng} there exist
$\mu_0<0$ and $f_0\in H^2_{per}([0,L])$ satisfying
$\mathcal{L}_1f_0=\mu_0f_0$ and so that
$\langle\mathcal{L}_1f_0,f_0\rangle<0$. Thus, by choosing $f=g=f_0$,
we obtain from \eqref{sa7},
$$
Q_R(f_0,f_0)=3\langle\mathcal{L}_1f_0,f_0\rangle<0.
$$

This implies that the first eigenvalue of  $\mathcal{A}_R$, say
$\lambda_1$, is negative. Now we will prove that the next eigenvalue
is the zero one. To do so, we will use the {\it min-max
characterization} of eigenvalues (see e.g., \cite[Theorem
XIII.1]{RSIV}). Thus, if $\lambda_2$ denotes the second eigenvalue
of $\mathcal{A}_R$, we have
\begin{equation}  \label{minmax2}
\lambda_2=\max_{(\phi_1,\phi_2)\in Y}\min_{(f,g)\in
Y\setminus\{(0,0)\}\atop{f\perp\phi_1, g\perp\phi_2}}
\frac{Q_R(f,g)}{\|(f,g)\|_{Y}^2}.
\end{equation}

By taking $\phi_1=\phi_2=f_0$, we see that
$$
\lambda_2\geq \min_{(f,g)\in Y\setminus\{(0,0)\}\atop{f\perp f_0,
g\perp f_0}} \frac{Q_R(f,g)}{\|(f,g)\|_{Y}^2}.
$$

Now, if $f\perp f_0$ and $g\perp f_0$ we obtain
$\langle\mathcal{L}_1f,f\rangle +
\langle\mathcal{L}_1g,g\rangle\geq0$ (recall that Lemma
\ref{lemmaAng} implies that $\mathcal{L}_1$ has a unique negative eigenvalue).
Moreover, since $\psi$ is a strictly positive function (and thus,
the last integral in \eqref{sa7} is non-negative) we obtain
$Q_R(f,g)\geq0$, which implies $\lambda_2\geq0$.

Finally, to prove that the third eigenvalue is strictly positive, we use the min-max
principle again, taking into account that $\mathcal{L}_1$ has a
unique negative eigenvalue and zero is a simple eigenvalue. This
proves part (i).\\

(ii) In this case, if $Q_I$ denotes the quadratic form associated
with $\mathcal{A}_I$, we have
\begin{equation} \label{sa8}
\begin{array}{lll}
Q_I(f,g)&:=&\displaystyle \langle\mathcal{A}_I(f,g),(f,g)\rangle\\\\
&=&\displaystyle\int_{0}^{L}\{ 2(-\partial_x^2+1-\psi)f^2+g^2 \}\,dx\\\\
&=& 2\langle\mathcal{L}_2f,f\rangle + \|g\|^2.
\end{array}
\end{equation}

Therefore, since $\mathcal{L}_2$ has no negative eigenvalue (see
Lemma \ref{lemmaAng}) we have $\langle\mathcal{L}_2f,f\rangle\geq0$
and then from \eqref{sa8} we deduce $Q_I(f,g)\geq0$. This implies
that $\mathcal{A}_I$ has no negative eigenvalue. Moreover, it is
easy to see from \eqref{s6} that zero is an eigenvalue with
eigenfunction $(\psi,0)$. This completes the proof of the theorem.\\
\fim

\subsection{Orbital stability}

In this subsection we prove our orbital stability result for the
periodic wave $(e^{it}\psi,\psi)$, where $\psi=\psi_1$ is the
cnoidal wave given in Theorem \ref{exist}. To make clear our notion
of orbital stability, we point out that  system \eqref{s1} has
translation and phase symmetries, i.e., if $(u(x,t),v(x,t))$ is a
solution for \eqref{s1}, so is
\begin{equation}\label{symmetries}
(e^{i\theta}u(x+x_0,t),v(x+x_0,t)),
\end{equation}
for any $\theta,x_0\in \mathbb{R}$. Thus, our notion of orbital
stability will be modulus such symmetries. To be more precise, we have the following definition

\begin{defi}    \label{orbitalstability}
A standing wave solutions for \eqref{s1} of the form \linebreak
$(e^{i\omega t}\psi_\omega(x),\phi_\omega(x))$, is said to be
orbitally stable in $X= H_{per}^1([0,L])\times
H_{per}^1([0,L])\times L_{per}^2([0,L])$ if for any $\varepsilon>0$
there exists $\delta>0$ such that if $(u_0,v_0,v_1)\in X$ satisfies
$||(u_0,v_0,v_1)-(\psi_\omega,\phi_\omega,0)||_X<\delta$, then the
solution $\overrightarrow{u}(t)=(u,v,v_t)$ of \eqref{s1} with
$\overrightarrow{u}(0)=(u_0,v_0,v_1)$ exists for all $t$ and
satisfies
$$
\displaystyle\sup_{t\geq0}\inf_{ s,y\in\mathbb{R}}
||\overrightarrow{u}(t)
-(e^{is}\psi_\omega(\cdot+y),\phi_\omega(\cdot+y),0)||_{X}<\varepsilon.
$$

Otherwise, $(e^{i\omega t}\psi_\omega(x),\phi_\omega(x))$ is said to
be orbitally unstable in $X$.
\end{defi}

From Theorem \ref{teoeigenR} we obtain the following properties

\begin{enumerate}
  \item[($i$)] The operator $\mathcal{A}$ has exactly one negative
  eigenvalue, that is, the nega\-tive eigenspace of $\mathcal{A}$, say
  $\mathcal{N}$, is one-dimensional.

  \item[($ii$)] For $\overrightarrow{f}=(\psi',\psi',0,0)$ and
  $\overrightarrow{g}=(0,0,\psi,0)$, the set
  $\mathcal{Z}=\{r_1\overrightarrow{f}+r_2\overrightarrow{g}; \;
  r_1,r_2\in \mathbb{R}\}$ is the kernel of operator $\mathcal{A}$.

  \item[($iii$)] There exists a closed subspace, say $\mathcal{P}$, such that
  $\langle\mathcal{A}u,u \rangle \geq \delta_0
\|u\|_X$ for all $u \in \mathcal{P}$ and some $\delta_0>0$.
\end{enumerate}

Therefore, from ($i$)--($iii$) we obtain the following orthogonal
decomposition for $X_{{\mathbb{R}}}=H_{per}^1([0,L])\times
H_{per}^1([0,L])\times H_{per}^1([0,L])\times L_{per}^2([0,L])$:
\begin{equation}  \label{5.1}
 X_{{\mathbb{R}}}=\mathcal{N} \oplus \mathcal{Z} \oplus \mathcal{P}.
\end{equation}

Next, for  $\omega \in \mathcal{I}=\left(4\pi^2/L^2, +\infty
\right)$ and $\psi_\omega$ the cnoidal wave given by Theorem
\ref{exist} we define $d:\mathcal{I}\rightarrow {\mathbb{R}}$ by
\begin{equation}    \label{5.2}
d(\omega)=\mathcal{E}(\Psi_\omega)+ \omega \mathcal{F}(\Psi_\omega),
\end{equation}
where, as before, $\Psi_\omega=(\psi_\omega,\psi_{\omega},0,0)$.

In the present setting, our orbital stability result in Theorem
\ref{stability} can be rephrased as follows

\begin{theo}
Let $\psi=\psi_1$ be the cnoidal wave given in Theorem \ref{exist}.
Then, the periodic traveling wave $(e^{it}\psi,\psi)$ is orbitally
stable in space $X$.
\end{theo}
\textbf{Proof.} Since the initial value problem associated with
\eqref{s1} is globally well-posed in $X$ (see Theorem \ref{t1.5}),
$\Psi_1=(\psi_1,\psi_1,0,0)$ satisfies \eqref{critical},
$X_{\mathbb{R}}$ admits the decomposition (\ref{5.1}) and
$\mathcal{N}$ is one-dimensional, the proof of the theorem follows
from the  {\it Abstract Stability Theorem} in Grillakis, Shatah and
Strauss \cite{GSS}, provided we are able to show that
$d''(\omega)>0$, where $d$ is the function defined in \eqref{5.2}.
 This was essentially proved in \cite{AL},
but for the sake of completeness we bring here the main steps. From
direct computation, we obtain $d'(\omega)=\mathcal{F}(\Psi_\omega)$.
Thus,
$$
d''(\omega)=\frac{d}{d\omega}\left(\int_0^L\psi_\omega^2(x)dx\right).
$$

But integrating \eqref{s6} over $[0,L]$, we get
$$
\int_0^L\psi_\omega^2(x)dx=\omega \int_0^L\psi_\omega(x)dx.
$$

Then, for the positivity of $d''(\omega)$ it  suffices  to show that
function
$$G(\omega)=\omega \int_0^L\psi_\omega(x)dx
$$
is strictly increasing.

In what follows we replace (up to a multiplicative positive
constant) $\eta$ with $\omega$ in the definition of $k$ and $\rho$
in Theorem \ref{exist}. Using that
$$\int_0^K cn^2(x;k)dx=\frac{[E(k)-(1-k^2)K(k)]}{k^2},
$$
(where $E(k)$ is the complete elliptic integral of the second kind)
$L=2\sqrt6K/\sqrt{\beta_3-\beta_1}$ and
$k^2=(\beta_3-\beta_2)/(\beta_3-\beta_1)$, we deduce
$$
\int_0^L\psi_\omega(x)dx=\beta_2L+24\frac{K}{L}[E-(1-k^2)K].
$$

Moreover, in view of the definitions of $k$ and $\rho$, we infer
that
$$
\beta_2=\frac{8K^2}{L}\left[\sqrt{k^4-k^2+1}+1-2k^2\right].
$$

As a consequence,
$$
\int_0^L\psi_\omega(x)dx=\frac{8K^2}{L}\left[\sqrt{k^4-k^2+1}-2+k^2\right]
+24\frac{KE}{L}\equiv H(k(\omega)).
$$

Finally,
$$
\frac{d}{d\omega}G(\omega)=\int_0^L\psi_\omega(x)dx+\omega
\frac{dH}{dk}\frac{dk}{d\omega}>0,
$$
where we have used that $k\mapsto H(k)$ is a strictly increasing
function and $dk/d\omega>0$ (see Theorem \ref{exist}). This
completes the proof of the theorem.\\
\fim

\subsection{Existence and Stability of non-explicit solutions}

In Subsection 6.1 we proved that system \eqref{s4} admits a periodic
wave solution for $\omega=1$ and $\psi_\omega=\phi_\omega$, where
$\psi_\omega$ is given explicitly by the formula in \eqref{s7}. The
advantage in that case is the reduction of system \eqref{s4} to a
single ordinary differential equation. However, one can naturally
ask if the system also admits a periodic solution for $\omega\neq1$.
In this regard, we shall prove that for $\omega$ sufficiently close
to $1$, system \eqref{s4} does admit an {\it even} periodic solution
such that at $\omega=1$ this solution is the aforementioned one. We
shall employ the Implicit Function Theorem combined with the spectral
results given in Theorem \ref{teoeigenR}.

Let $H^s_{per,e}([0,L])$ be the subspace of
$H^s_{per}([0,L])$ constituted by the even distributions and denote
$X_e=H^2_{per,e}([0,L])\times H^2_{per,e}([0,L])$ and
$Y_e=L^2_{per,e}([0,L])\times L^2_{per,e}([0,L])$. Define the
function $\Phi:\mathbb{R}\times X_e\rightarrow Y_e$ by
$$
\Phi(\omega,\psi,\phi)=(-\psi''+\omega\psi-\psi\phi,
-\phi''+\phi-\psi^2).
$$

In view of Theorem \ref{exist} we deduce that
$\Phi(1,\psi_1,\psi_1)=(0,0)$. Moreover, if $\Phi_{(\psi,\phi)}$
denotes the Fr\'echet derivative of $\Phi$ at $(\psi,\phi)$, it is
easy to check that
\begin{equation*}
\displaystyle\Phi_{(\psi,\phi)}(\omega,\psi,\phi)=\left(
\begin{array}{cccc}
\displaystyle -\partial_x^2+\omega-\phi & -\psi\\\\
-2\psi & \displaystyle -\partial_x^2+1
\end{array}\right).
\end{equation*}

Thus, at $\omega=1$ and $\psi=\phi=\psi_1$, we obtain
\begin{equation*}
\displaystyle\mathcal{B}:=\Phi_{(\psi,\phi)}(1,\psi_1,\psi_1)=\left(
\begin{array}{cccc}
\displaystyle -\partial_x^2+1-\psi_1 & -\psi_1\\\\
-2\psi_1 & \displaystyle -\partial_x^2+1
\end{array}\right).
\end{equation*}

Let us prove that $\mathcal{B}$ is a bijection from $X_e$ into
$Y_e$. In fact, it is sufficient to show that 0 does not belong to
$\sigma(\mathcal{B})$. An elementary calculation shows us that
$(f,g)\in \mbox{Ker}(\mathcal{B})$ if and only if $(f,g)\in
\mbox{Ker}(\mathcal{A}_R)$, where $\mathcal{A}_R$ is the operator
given by \eqref{sa5}. But, from Theorem \ref{teoeigenR} we have
$\mbox{Ker}(\mathcal{A}_R)=[(\psi_1',\psi_1')]$ (as an operator on
$L^2_{per}([0,L])\times L^2_{per}([0,L])$). However, since $\psi_1$
is an even function, it follows that $\psi_1'\not\in
L^2_{per,e}([0,L])$ and so $0\not \in \sigma(\mathcal{B})$ (as an
operator on $Y_e$).

Consequently, from the Implicit Function Theorem there exist an
$\varepsilon>0$ and a unique smooth function
$\digamma:(1-\varepsilon,1+\varepsilon)\rightarrow X_e$,
$$
\digamma(\omega)=(\psi_\omega,\phi_\omega),
$$
such that $\digamma(1)=(\psi_1,\psi_1)$ and
$\Phi(\omega,\digamma(\omega))=(0,0)$, for all $\omega\in
(1-\varepsilon,1+\varepsilon)$, that is, the pair
$(\psi_\omega,\phi_\omega)$ is a solution of the system \eqref{s4}.

\begin{rema}
The periodic solution we found here are also orbitally stable. This
can be proved by using classical perturbation theory (see
\cite{KATO}) to show that the linearized operators arising in this
context have the same spectral properties as those ones in Theorem
\ref{teoeigenR} (for related references see e.g. \cite{AMP1},
\cite{NP1} and references therein).
\end{rema}


\centerline{\textbf{Acknowledgment}}

The authors want to take the opportunity to thank Professor F\'abio Natali for 
helpful discussion concerning this work.



E-mail: apastor@impa and farah@impa.br
\end{document}